\documentclass[11pt,psfig]{article}
\usepackage{amsmath,amssymb,epsfig, epstopdf, color}
\usepackage{subfigure}

\newtheorem{theorem}{Theorem}[section]
\newtheorem{definition}[theorem]{Definition}
\newtheorem{proposition}[theorem]{Proposition}
\newtheorem{lemma}[theorem]{Lemma}

\newtheorem{remark}[theorem]{Remark}

\newtheorem{assumption}[theorem]{Assumption}

\newcommand{\g}{\gamma}
\newcommand{\e}{\epsilon}
\newcommand{\aone}{\alpha_1}
\newcommand{\atwo}{\alpha_2}

\newcommand{\lam}{\lambda}

\makeatletter
\@addtoreset{equation}{section}

\makeatother
\begin{document}

\begin{center}
{\Huge \textbf{Time consistent portfolio management }}
\vspace{-0.2cm}

\mbox{}\\[0pt]
 Ivar Ekeland\\[0pt]
CEREMADE et Institut de Finance\\[0pt]
Universite Paris-Dauphine\\[0pt]
75775 Paris CEDEX 16\\[0pt]
ekeland@math.ubc.ca\\[0pt]

\vspace{0.2cm} Oumar Mbodji \\[0pt]
Department of Mathematics \& Statistics\\[0pt]
McMaster University \\[0pt]
1280 Main Street West \\[0pt]
Hamilton, ON, L8S 4K1\\[0pt]
oumarms@math.mcmaster.ca

\vspace{0.2cm} Traian A.~Pirvu  \footnote{%
Work supported under NSERC grant 298427-04. We thank the referees for valuable advice, suggestions and a thorough reading of the first version.} \\[0pt]
Department of Mathematics \& Statistics\\[0pt]
McMaster University \\[0pt]
1280 Main Street West \\[0pt]
Hamilton, ON, L8S 4K1\\[0pt]
tpirvu@math.mcmaster.ca \vspace{1cm}

\mbox{}\\[0pt]

\end{center}
\vspace{-1.4cm}
\noindent \textbf{Abstract.} This paper considers the portfolio
management problem for an investor with finite time horizon who is allowed to consume and take out life insurance. Natural assumptions, such as different discount rates for consumption and life insurance lead to time inconsistency. This situation can also arise when the investor is in fact a group, the members of which have different utilities and/or different discount rates. As a consequence,
the optimal strategies are not implementable. We focus on hyperbolic discounting,
which has received much attention lately, especially in the area of
behavioural finance. Following \cite{EkePir}, we consider the resulting problem as a leader-follower
game between successive selves, each of whom can commit for an infinitesimally small amount of time.
 We then define policies as subgame perfect equilibrium strategies. Policies are characterized by an
integral equation which is shown to have a solution in the case of CRRA utilities. 
Our results can be extended for more
general preferences as long as the equations admit solutions.
 Numerical simulations reveal that for the Merton problem with hyperbolic discounting, the
consumption increases up to a certain time, after which it decreases; this
pattern does not occur in the case of exponential discounting, and is
therefore known in the litterature as the ``consumption puzzle". Other numerical experiments explore the
effect of time varying aggregation rate on the insurance premium.

\vspace{0.4cm}



\section{Introduction}

The investment/consumption problem in a stochastic context was considered by
Merton \cite{Mer69} and \cite{Mer71}. His model consists in a risk-free
asset with constant rate of return and one or more stocks, the prices of
which are driven by geometric Brownian motion. The horizon $T$ is
prescribed, the portfolio is self-financing, and the investor seeks to
maximize the expected utility of intertemporal consumption plus the final
wealth. Merton provided a closed form solution when the utilities are of
constant relative risk aversion (CRRA) or constant absolute risk aversion
(CARA) type. It turns out that for (CRRA) utilities the fraction of wealth
invested in the risky asset is constant through time. Moreover for the case of (CARA) utilities, they are
linear in wealth.

Richard \cite{Rich} added life insurance to the investor's portfolio, assuming
an arbitrary but known distribution of death time. In the same vein Pliska
\cite{Pliska} studied optimal life insurance and consumption for an income
earner whose lifetime is random and unbounded. More recently Kwak et al.\cite%
{Kwak} looked at the problem of finding optimal investment and consumption
for a family whose parents receive deterministic labor income until some
deterministic time horizon.

The aim of this paper is to revisit these problems in the case when the
psychological discount rate is not constant. By now there is substantial evidence that people discount the future at a
non-constant rate. More precisely, there is experimental evidence (see
Frederick et. al. \cite{Frederick} for a review) that people are more sensitive
to a given time delay if it occurs earlier: for instance, a person might prefer to
get two oranges in 21 days than one orange in 20 days, but also prefer to
get one orange right now than two oranges tomorrow. This is known as \textbf{%
the common difference effect}, and would not occur if future
utilities are discounted at a constant rate. Individual behaviour is best described by
\textbf{hyperbolic discounting}, where the discount factor is $h(t)=(1+at)^{-%
\frac{b}{a}},$ with $a,b>0$. The corresponding discount rate is $r\left(
t\right) =\frac{b}{ 1+at} $, which starts from $r\left( 0\right) =\frac{b%
}{a}$ and decreases to zero. Because of its empirical support, hyperbolic
discounting has received a lot of attention in the areas of: microeconomics,
macroeconomics and behavioural finance. We just mention here among others
the works of Loewenstein and Prelec \cite{LoPre}, Laibson \cite{Lai} and
Barro \cite{Bar}.

It is well-known that, for non-constant discount rates, optimal strategies are time inconsistent: for $t_1 < t_2$, the planner at time $t_1 $ will find a strategy $f_1$ to be optimal on $[t_1, \infty)$, while the planner at time $t_2$ will find a different strategy $f_2$ to be optimal on that interval. As a result, the planner at time $t_2$ will not implement the strategy devised by the planner at time $t_1$, unless there exists some commitment mechanism. If there is none, then the strategy $f_1$, which is optimal from the perspective of the planner at time $t_1$, is not implementable, and the planner at time $t_1$ must look for a second-best strategy. This situation was first analyzed by Strotz \cite{STRO}, and
this line of research has been pursued by many others (see Pollak \cite{Pol}%
, Phelps \cite{Phelps}, Peleg and Yaari \cite{PeYa}, Goldmann \cite{Gol},
Laibson \cite{Lai}, Barro \cite{Bar}, Krusell and Smith \cite{KruSm}),
mostly in the framework of planning a discrete-time economy with production
(Ramsey's problem). It is by now well established that time-consistent
strategies are Stackelberg equilibria of a leader-follower game among
successive selves (today's self has divergent interests from tomorrow's).\
More recently, the problem has been taken up again by Karp \cite{Karp2},
\cite{Karp}, \cite{Karp-Fuji}, \cite{Karp-Lee}. Luttmer and Mariotti, \cite%
{LutMar}, Ekeland and Lazrak \cite{EkeLaz}, \cite{Ekl1}, \cite{EkL}, always
within the framework of planning economic growth. In a series of papers
Bj\"{o}rk and Murgoci \cite{Bj1}, Bj\"{o}rk, Murgoci and Zhou \cite{Bj2} look at
the mean variance problem which is also time inconsistent.

Ekeland and Pirvu \cite{EkePir} seem to have been the first to have
considered the Merton problem with non-constant psychological discount rates. They studied
the case of an investor who has a CRRA utility $u\left( c\right) =-\frac{1}{p%
}c^{p},\ p<1$, and a quasi-exponential discount factor $h\left( t\right) $,
that is, $h\left( t\right) $ must belong to one of the families:%
\begin{eqnarray*}
h\left( t\right) &=&\lambda e^{-r_{1}t}+\left( 1-\lambda \right) e^{-r_{2}t},
\\
h\left( t\right) &=&\left( 1+at\right) e^{-rt}.
\end{eqnarray*}
Extending the basic idea of Ekeland and Lazrak \cite{EkeLaz} to the
stochastic framework, they find time-consistent strategies in the limiting
case when the investor can commit only during an infinitesimal time
interval. They show that time-consistent strategies exist if a certain BSDE
has a solution, and they show that, because of the special form of the
discount factor, this BSDE reduces to a system of two ODEs which has a
solution.

The aim of this paper is to extend these results to more general discount
rates, and more general problems. Quasi-exponential discount rates, although
mathematically convenient, are not realistic. As we saw earlier, empirically
observed discount rates among individuals tend to be hyperbolic. But there
is another, perhaps more compelling, reason why general discount rates are
of interest. Standard portfolio theory assumes that the investor is an
individual. However, in most situations investment decisions are made by a
group, such as the management team in the case when the investor entrusts
his portfolio to professionals. Even when the investor manages the portfolio
directly, the word "investor" which is suggestive of a single decision-maker
very often hides a different reality, namely the family: one would expect
the husband and the wife to take part in investment decisions concerning the
couple. By now, the relevant economic literature has made abundantly clear
(see Chiappori and Ekeland, \cite{Chia}) that the group cannot be represented
by a single utility function. Instead, there should be one utility function,
and one discount factor per member of the group. The actual decision taken
is the result of negociations within the group, a kind of black box which
cannot be opened by outsiders. However, if the group is efficient, that is,
if the outcome is Pareto optimal, then it can be modelled by maximising a
suitable convex combination of the members' utilities, the weight conferred
to each individual representing his/her power within the group. In the (very) particular case when all members of the group have the same utility, but different discount rates, the group behaves as a single individual with non-constant discount rate.

The difficulty in dealing with non-constant discount rates is to define time-consistent strategies and to
prove that they exist. We follow the approach pionneered by Ekeland and
Lazrak \cite{EkeLaz} in the deterministic framework, by considering the
limiting case when the decision-maker can commit only during an
infinitesimal amount of time. This approach was already followed in our
earlier work \cite{EkePir}, in the case of quasi-exponential discount
factors. The proofs in that paper do not readily extend to the case of
general discount factors, so in the present work we present a different
method. Whereas \cite{EkePir} characterizes the time-consistent stategies in
terms of a certain BSDE, we now characterize them in terms of a certain
"value function", which is shown to satisfy a certain integral equation
which has a natural interpretation. Assuming utilities to be (CRRA), we
decouple time and space, and reduce it to a one-dimensional integral
equation, which we solve by a fixed-point argument. Moreover this one
dimensional equation is amenable to numerical treatments so one can compare
the equilibrium policies arising from different choices of discounting. The numerical scheme we employ
consists of the discretization of the one dimensional equation in three steps. This is based on a Riemann
sum approximation of the integral. We obtain closed form solutions in certain cases.

We show that hyperbolic discounting may result in consumption patterns which
are observationally different from the optimal strategies in the Merton
model. The latter predicts that consumption grows smoothly over time if the
interest rate exceeds the discount rate (or decays smoothly otherwise).
However, household data indicates that consumption is hump-shaped. This is
referred to in the literature as the consumption puzzle, and we show that it
can arise as a time-consistent strategy in certain cases of hyperbolic
discounting.

By running numerical simulations we study the effect of the weight given by the insurer to the beneficiaries on the life insurance process.

\textbf{Organization of the paper}: The remainder of this paper is organized
as follows. In section $2$ we describe the model and formulate the
objective. Section $3$ introduces the value function. Section 4 presents
the main result. Section $5$ deals with CRRA utilities. Numerical results are discussed in Section $6$. An extension to multiple managers is discussed in Section $7.$ The paper ends with an appendix containing some proofs.

\section{The Model}

\subsection{The decisions}

Consider a financial market consisting of a savings account and one stock
(the risky asset). The inclusion of more risky assets can be achieved by
notational changes. The savings account accrues interest at the riskless
rate $r>0.$ The stock price per share follows an exponential Brownian motion
\begin{equation*}
dS(t)=S(t)\left[ \alpha \,dt+\sigma \,dW(t)\right] ,\quad 0\leq t\leq \infty
,
\end{equation*}%
where $\{W(t)\}_{t\geq 0}$ is a $1-$dimensional Brownian motion on a
filtered probability space,\\ $(\Omega ,\{\mathcal{F}_{t}\}_{t\geq 0},\mathbb{P}).$ The filtration $\{\mathcal{F}_{t}\}$ is the completed
filtration generated by $\{W(t)\}$. Let us denote by $\mu \triangleq \alpha
-r>0$ \emph{the excess return}.

A decision-maker in this market is continuously investing in the stock and
the bond, consuming and buying life insurance, while receiving income at the
continuous deterministic rate $i(t).$ This assumption is key in deriving our results.
Relaxing it to accommodate for problems relevant to small enterprises  is not obvious
and would be an interesting research project.

 Life insurance is offered as a
sucession of term contracts with infinitesimally small horizon. At every
time $t$, a contract is offered, costing $1$ unit of account. If the holder
dies immediately after, the insurance company pays $l\left( t\right) $ to
his/her beneficiaries.\ The deterministic function $l\left( t\right) $ is
prescribed.

At every time $t$, the investor chooses $\zeta (t)$, the investment in the
risky asset,   $c(t)$ the consumption,  and $p(t)$, the amount of life
insurance. Given an adapted process $\{\zeta (t),c(t),p(t)\}_{t\geq 0}$, the
equation describing the dynamics of wealth ${X^{\zeta ,c,p}(t)}$ is given by
\begin{eqnarray}
dX^{\zeta ,c,p}(t) &=&rX^{\zeta ,c,p}(t)dt-c(t)dt-p(t)dt+i(t)dt+\zeta
(t)(\alpha \,dt+\sigma dW(t))  \notag  \label{equ:wealth-one} \\
X^{\zeta ,c,p}(0) &=&X\left( 0\right),
\end{eqnarray}%
the initial wealth $X(0)$ being exogenously specified.

We assume a benchmark deterministic time horizon $T$. The investor is
alive at time $t=0$ and has a lifetime denoted by $\tau $, which is a
non-negative random variable defined on the probability space $(\Omega ,%
\mathcal{F},\mathbb{P})$ and independent of the Brownian motion $W.$ Denote
by $g\left( t\right) $ its density and by $G\left( t\right) $ its
distribution:

\begin{equation*}
G(t)\triangleq\mathbb{P}(\tau <t)=\int_{0}^{t}g(u)\,du.
\end{equation*}

It will be useful for later computations to introduce the hazard function $%
\lambda \left( t\right) $, that is, the instantaneous death rate, defined by

\begin{equation*}
\lambda (t)\triangleq\lim_{\delta t\rightarrow 0}\frac{\mathbb{P}(t\leq \tau
<t+\varepsilon \ |\ \tau \ \geq t)}{\varepsilon }=\frac{g(t)}{1-G(t)},
\end{equation*}%
so that $g(t)=\lambda (t)\exp \{-\int_{0}^{t}\lambda (u)\,du\}$. We have,
from the definition:

\begin{equation}\label{*0}
\mathbb{P}(\tau <s\ |\ \tau >t)=1-\exp \{-\int_{t}^{s}\lambda (u)\,du\}.
\end{equation}

and

 \begin{equation}\label{*1}
 \mathbb{P}(\tau>T | \tau>t)= \exp\{-\int_{t}^{T} \lambda(u)\,du \}.
 \end{equation}

Next we turn to risk preferences.

\subsection{Utility functions}

\begin{definition}\label{util}
 A utility function $U$ is a strictly increasing, strictly concave differentiable real-valued function
defined on $[0,\infty)$ which satisfies the Inada conditions

\begin{equation}\label{in}
U'(0)\triangleq \lim_{x\downarrow 0} U'(x)=\infty,\qquad U'(\infty)\triangleq \lim_{x\rightarrow \infty} U'(x)=0.
\end{equation}

\end{definition}

The strictly decreasing $C^1$ function $U'$ maps $(0,\infty)$ onto $(0,\infty)$ and hence has a strictly
decreasing, $C^{1}$ inverse $I: (0,\infty)\rightarrow (0,\infty).$

The legacy process of the decision-maker, $\{Z^{\zeta ,c,p}(t)\}_{t\geq 0},$ is defined by
\begin{equation}\label{le}
Z^{\zeta ,c,p}(t)\triangleq\eta (t)X^{\zeta ,c,p}(t)+l(t)p(t),
\end{equation}
 where $\eta \left( t\right) $ and $%
l\left( t\right) $ \ are prescribed deterministic and continuous functions.
The legacy is the sum of two terms:\ the first one, $\eta (t)X^{\zeta ,c,p}(t)$, is
the part of his wealth which will benefit his heirs (after taxes, and
various costs), and the second one $l(t)p(t)$ is the life insurance.
Although the insurance premium $p(t)$ is allowed to be negative we require
that the legacy $Z^{\zeta ,c,p}(t)$ stays positive. A negative $p(t)$ means that the decision maker can sell life insurance.

Let $U_1, U_2, U_3$ be utility functions as in Definition \ref{util}; $U_1$ is the utility from intertemporal consumption,
$U_2$ is the utility of the final wealth and $U_3$ is the utility of the legacy.

Next, we define the admissible strategies. Sometimes, to ease notations, we write ${X}^{t,x}(s)$ and ${Z}^{t,x}(s)$ for $%
\mathbb{E}[{X}^{\zeta ,c,p}(s)|{X}^{\zeta ,c,p}(t)=x]$ and $\mathbb{E}[{Z}%
^{\zeta ,c,p}(s)|{X}^{\zeta ,c,p}(t)=x].$

\begin{definition} \label{def:portfolio-proportions} An $\mathbb{R}^{3}$-valued
stochastic process $\{\zeta (t),c(t),p(t)\}_{t\geq 0}$ is called an
admissible strategy process if

\begin{itemize}
\item it is progressively measurable with respect to the sigma algebra $%
\sigma (\{\ W(t)\}_{t\geq 0})$,

\item $c(t)\geq 0,Z^{\zeta ,c,p}(t)\geq 0\,\,\mbox{for all}\,,\text{ a.s.};$
${X}^{\zeta ,c,p}(T)\geq 0,\text{ a.s.}$

\item moreover we require that for all $t, x\geq 0$
\begin{equation}  \label{189}
\mathbb{E} \sup_{\{t\leq s\leq T\}} |U_{1}(c(s))|<\infty,\,\,\mathbb{E}
|U_{2}({X}^{t,x}(T))|  <\infty,\,\, \mathbb{E}\sup_{\{t\leq s\leq T\}} |U_{3}(Z^{t,x}(s))|<\infty.
\end{equation}
\end{itemize}
\end{definition}

The last set of inequalities are purely technical and are satisfied for e.g.
bounded strategies. They are essential in proving our main result and are related to the fact that the expected utility criterion is continuously updated. \subsection{The intertemporal utility}

In order to evaluate the performance of an investment-consumption-insurance
strategy the decision maker uses an expected utility criterion. For an
admissible strategy process\\ $\{{\zeta }(s),{c}(s),p(s)\}_{s\geq 0}$ and its
corresponding wealth process $\{X^{\zeta ,c,p}(s)\}_{s\geq 0},$ we denote the expected
intertemporal utility  by
\begin{eqnarray}
J(t,x,\zeta ,c,p) &\triangleq &\mathbb{E}\bigg[\int_{t}^{T\wedge \tau
}h(s-t)U_{1}(c(s))\,ds+nh(T-t){U}_{2}(X^{\zeta
,c,p}(T))1_{\{\tau >T|\tau >t\}}  \notag \\
&+&m(\tau -t)\hat{h}(\tau -t){U}_{3}(Z^{\zeta ,c,p}(\tau ))1_{\{\tau
\leq T|\tau >t\}}\bigg|X^{\zeta ,c,p}(t)=x\bigg],  \label{01FUNCT}
\end{eqnarray}%
where:

\begin{itemize}

\item $n>0$ is a constant

\item $m\left( t\right) >0$ is a continuous function

\item $h$ and $\hat{h}$ are continuously differentiable, positive and
decreasing functions, such that $h\left( 0\right) =\hat{h}\left( 0\right) =1.$%

\end{itemize}

The interpretation is as follows. The decision-maker will collect $X^{\zeta ,c,p}(T)$
at time $T$, if he is still alive at time $T$, and the coefficient $n$ is
the weight he attributes to getting that lump sum, as compared to the utility of
continuous consumption up to time $T$. The function $h\left( t\right) $ is
his discount function, and it is no longer restricted to the exponential and
quasi-exponential type.

He may, however, die before time $T$, in which case his wealth will accrue
to others, and the decision-maker is taking the utility of his beneficiaries
into account when managing his portfolio.

Since the death time $\tau $ is independent of the uncertainty driving the
stock, we have the following simplified expression for the functional $J,$
which is proved in Appendix A.

\begin{lemma}
\label{L1} \label{ain200} The functional $J$ of \eqref{01FUNCT} equals
\begin{eqnarray}  \label{w9}
J(t,x,\zeta ,c,p)&=& \mathbb{E}\bigg[ \int_{t}^{T} Q(s,t)
U_{1}(c(s))\,ds \\
&+&\int_{t}^{T} q(s,t){U}_{3}(Z^{t,x}(s))\,ds+n Q(T,t) {U}_{2}(X^{t,x}(T))\bigg],
\end{eqnarray}

where%
\begin{eqnarray}
q(s,t) &\triangleq &\bar{h}(s-t)\lambda (s)\exp \{-\int_{t}^{s}\lambda
(z)\,dz\},\quad \bar{h}(t)\triangleq m(t)\hat{h}(t)  \label{((o} \\
Q(s,t) &\triangleq &h(s-t)\exp \{-\int_{t}^{s}\lambda (z)\,dz\}  \label{((o0}
\end{eqnarray}
\end{lemma}

A natural objective for the decision maker is to maximize the above expected
utility criterion. However, because neither $q$ nor $Q$ are exponentials,
time inconsistency will bite, that is, a strategy that will be considered to
be optimal at time $0$ will not be considered so at later times, so it will
not be implemented unless the decision-maker at time $0$ can constrain the
decision-maker at all times $t>0$.

\subsection{Time-consistent strategies}

We now introduce a special class of time-consistent strategies, which will
henceforth be called \emph{policies}. That is, we consider that the
decision-maker at time $t$ can commit his successors up to time $\varepsilon
$, with $\varepsilon \rightarrow 0$, and we seek strategies which it is
optimal to implement right now conditioned on them being implemented in the
future. 

More precisely, suppose that a strategy $f$ is time-consistent. This means that, if it has been applied up to time $t$, the decision-maker at time $t$ will apply it as well. Since there is no commitment mechanism to force him to do so, he will only apply strategy $f$ if it is in his own best interests. Denote his current wealth by $X(t)$. He has two possibilities: either to stick to the strategy $f$, or to apply another one. To simplify matters, we will assume that the decision-maker considers only a very short time interval, $[t, t+\epsilon]$, so short in fact that all strategies can be assumed to be constant on that interval. The decision-maker then just compares the effect of investing $\bar{\zeta}$, consuming $\bar{c}$ and buying $\bar{p}$ worth of insurance, as required by the strategy $f$ at time t, with the effect of investing $\zeta$, consuming $c$ and buying $p$ worth of insurance, for different (constant) values. There will be, as usual, an immediate effect, corresponding to the change in consumption between $t$ and $t+\epsilon$, and a long-term effect, corresponding to the change in wealth at time $t+\epsilon$.

Let us formalize this idea:

\begin{definition}
\label{finiteh}An admissible trading strategy $\{\bar{\zeta}(s),\bar{c}(s),%
\bar{p}(s)\}_{t\leq s\leq T}$ is a \emph{policy} if there exists a map $%
F=(F_{1},F_{2},F_{3}):[0,T]\times \mathbb{R}\rightarrow \mathbb{R}\times
\lbrack 0,\infty )\times \mathbb{R}$ such that for any $t,x>0$
\begin{equation}
{\lim \inf_{\epsilon \downarrow 0}}\frac{J(t,x,\bar{\zeta},\bar{c},\bar{p}%
)-J(t,x,\zeta _{\epsilon },c_{\epsilon },p_{\epsilon })}{\epsilon }\geq 0,
\label{opt}
\end{equation}%
where:%
\begin{equation}\label{0000eq}
\bar{\zeta}(s)={F_{1}(s,\bar{X}(s))},\quad \bar{c}(s)=F_{2}(s,\bar{X}%
(s)),\quad \bar{p}(s)=F_{3}(s,\bar{X}(s))
\end{equation}%
and the wealth process $\{\bar{X}(s)\}_{s\in \lbrack t,T]}$ is a solution of
the stochastic differential equation (SDE):
\begin{equation}
d\bar{X}(s)=[r\bar{X}(s)+\mu F_{1}(s,\bar{X}(s))-F_{2}(s,\bar{X}(s))-F_{3}(s,%
\bar{X}(s))+i(s)]ds+\sigma F_{1}(s,\bar{X}(s))dW(s).  \label{0dyn}
\end{equation}

Here, the process $\{{\zeta }_{\epsilon }(s),{c}_{\epsilon
}(s),p_{\epsilon }(s)\}_{s\in \lbrack t,T]}$ is another
investment-consumption strategy defined by
\begin{equation}
\zeta _{\epsilon }(s)=%
\begin{cases}
\bar{\zeta}(s),\quad s\in \lbrack t,T]\backslash E_{\epsilon ,t} \\
\zeta (s),\quad s\in E_{\epsilon ,t},%
\end{cases}
\label{1e}
\end{equation}%
\begin{equation}
c_{\epsilon }(s)=%
\begin{cases}
\bar{c}(s),\quad s\in \lbrack t,T]\backslash E_{\epsilon ,t} \\
c(s),\quad s\in E_{\epsilon ,t},%
\end{cases}
\label{2e}
\end{equation}%
\begin{equation}
p_{\epsilon }(s)=%
\begin{cases}
\bar{p}(s),\quad s\in \lbrack t,T]\backslash E_{\epsilon ,t} \\
p(s),\quad s\in E_{\epsilon ,t},%
\end{cases}
\label{3e}
\end{equation}
with $E_{\epsilon,t}=[t,t+\epsilon],$ and $\{{\zeta}(s),{c}(s),{p}%
(s)\}_{s\in E_{\epsilon,t} }$ is any strategy for which $\{{\zeta}%
_{\epsilon}(s),{c}_{\epsilon} (s), {p}_{\epsilon} (s)\}_{s\in[t,T]}$ is an
admissible policy.
\end{definition}

In other words, policies are Markov strategies such that unilateral
deviations during an infinitesimally small time interval are penalized. Note
that:

\begin{itemize}
\item this does not mean that unilateral deviations during a finite interval
of time are penalized as well: it is possible that deviating from the policy
between $t_{1}$ and $t_{2}$ will be to the advantage of all the
decision-makers operating between $t_{1}$ and $t_{2}.$

\item however, if a Markov strategy is not a policy, then it certainly will
not be implemented, for at some point it will be to the advantage of some lone decision-maker to deviate,
during a very small time interval, which is enough to compromise all the plans laid by his predecessors.
\end{itemize}

So time-consistency in the sense of Definition \ref{finiteh} is a minimal
requirement for rationality: policies are the only Markov strategies that
the decision-maker should consider.

\section{The Value Function}

We now extend to this situation the notion of a value function, which is classical in optimal control.
 Let $m\triangleq m(0),$ and $I_{1},  I_{3}$ be the inverse functions of $U'_{1}, U'_{3}.$

\begin{definition}
\label{de1}
Let $v:[0,T]\times \mathbb{R}\rightarrow \mathbb{R}$ be a $%
C^{1,2}$ function, concave in the second variable. We shall say that $v$ is
a value function if we have:%
\begin{equation}
v(t,x)=J(t,x,\bar{\zeta},\bar{c},\bar{p}).  \label{100ie19(*}
\end{equation}%
Here the admissible process $\{\bar{\zeta}(s),\bar{c}(s),\bar{p}(s)\}_{s\in
\lbrack t,T]}$ is given by: 
\begin{equation}
\bar{\zeta}(s)\!\triangleq\! {-\frac{\mu \frac{\partial v}{\partial x}(s,\bar{X}(s))}{%
\sigma ^{2}\frac{\partial ^{2}v}{\partial x^{2}}(s,\bar{X}(s))}},\,\, \bar{c}%
(s)\!\triangleq\! I_{1}\!\left( \frac{\partial v}{\partial x}(s,\bar{X}(s))\right) ,\,\, 
\bar{p}(s)\!\triangleq\! \frac{1}{l(s)}\!\left[ I_{3}\left( \frac{1}{m}\frac{%
\partial v}{\partial x}(s,\bar{X}(s))\right)\! -\!\eta (s)\bar{X}(s)\right]\!,  \label{ie}
\end{equation}%
where $\bar{X}(s)$ is the corresponding wealth process defined by the SDE 
\begin{eqnarray}  \label{sde0}
\!\!\!\!\!\!d\bar{X}(s)&=&\bigg[r\bar{X}(s)-\frac{%
\mu ^{2}\frac{\partial v}{\partial x}(s,\bar{X}(s))}{\sigma ^{2}\frac{%
\partial ^{2}v}{\partial x^{2}}(s,\bar{X}(s))}-I_{1}\left( \frac{%
\partial v}{\partial x}(s,\bar{X}(s))\right) \!\!   \\\notag &-&   \!\!\frac{1}{l(s)}\left[
I_{3}\left( \frac{1}{m}\frac{\partial v}{\partial x}(s,\bar{X}(s))\right)
-\eta (s)\bar{X}(s)\right] +i(s)\bigg]ds-\frac{\mu \frac{\partial v}{%
\partial x}(s,\bar{X}(s))}{\sigma \frac{\partial ^{2}v}{\partial x^{2}}(s,%
\bar{X}(s))}dW(s)  \label{1o*} \\
\bar{X}\left( t\right) \!\!\!\!\! &=&\!\!\!\!\!x.
\end{eqnarray}
\end{definition}
The economic interpretation is very natural: if one applies the Markov
strategy associated with $v$ by the relations (\ref{ie}), and computes the
corresponding value of the investor's criterion starting from $x$ at time $t$%
, one gets precisely $v\left( t,x\right) $. In other words this is
fundamentally a fixed-point characterization.

 Let us define the functions $F_{1},F_{2},F_{3}$ by:%
\begin{equation}
\!\!\!\!\!F_{1}(t,x)\triangleq-\frac{\mu \frac{\partial v}{\partial x}(t,x)}{\sigma
^{2}\frac{\partial ^{2}v}{\partial x^{2}}(t,x)},\,F_{2}(t,x)\triangleq I_{1} \left( \frac{%
\partial v}{\partial x}(t,x)\right),\,  F_{3}(t,x)\triangleq%
\frac{1}{l(t)}\left[ I_{3} \left( \frac{1}{m}\frac{\partial v}{\partial x}%
(t,x)\right)-\eta (t)x\right].  \label{109con}
\end{equation}%
Next we impose a technical assumption; for a $C^{1,2}$ function
$f:[0,T]\times \mathbb{R}\rightarrow \mathbb{R},$ let us define
the operator $L f$ by
$$ L f(t,x)\triangleq \frac{\partial f}{\partial t}(t,x)+(r x+\mu F_{1}(t,x)-F_{2}(t,x)-F_{3}(t,x)+i(t))\frac{\partial f}{\partial x}(t,x)+\frac{\sigma^{2} F_{1}^{2}(t,x)}{2} \frac{\partial^2 f}{\partial x^2}(t,x).$$

\begin{assumption}\label{A2}
Assume that the PDEs
 \begin{equation}\label{09}
 L f(t,x)=0,\quad f(s,x)=g(x),
 \end{equation}
 have a $C^{1,2}$ solution on $[0,s]\times \mathbb{R}\rightarrow \mathbb{R}$ with exponential growth. Here $t<s\leq T,$ and $g(x)$ is
 one of the functions
 $$ U_{1}(F_{2}(s,x)):\,\, t<s\leq T,\quad U_{3}(\eta (s)x+l(s)F_{2}(s,x)):\,\, t<s\leq T,\quad  U_{2}(x):\,\, s=T.$$
 \end{assumption}


\section{Main Result}

The following Theorem states the central result of our paper. It involves the notions of policies and value
function for which we gave economic intuition.

\begin{theorem}\label{existence}
Let $v$ be a value function which satisfies Assumption \ref{A2}. Then, $\{\bar{\zeta}(s),\bar{c}(s),\bar{p}(s)\}_{s\in
\lbrack t,T]}$ given by \eqref{ie} is a policy.
\end{theorem}

We proceed in two steps. First we show that the value function $v$ satisfies a partial differential equation with a non-local term and this is done in the following Lemma, which is proved in Appendix B.

\begin{lemma}\label{pDe}
 The function $v$ solves the following equation
 \begin{equation}  \label{12dE}
\!\!\!\!\!\!\!\!\frac{\partial {v}}{\partial t}(t,x)+\bigg(rx+ \mu F_1(t,x) -F_2(t,x)-F_3(t,x) +i(t)  \bigg) \frac{\partial {v}}{\partial x}(t,x)+
\end{equation}
$$\frac{\sigma^{2} F^{2}_{1}(t,x)} {2}\frac{\partial^{2} v}{\partial x^{2}}(t,x)+{U}_{1}(F_2(t,x))+m{U}_{3}(x+l(t)F_3(t,x))= $$

$$
 \mathbb{E}\bigg[
\int_{t}^{T}  \frac{\partial Q}{\partial{t}} (s,t) U_{1}(F_2(s, \bar{X}^{t,x}(s) )  )\,ds+\int_{t}^{T}\frac{\partial q}{\partial{t}} (s,t){U}_{3}(\bar{Z}^{t,x}(s))\,ds+n \frac{\partial Q}{\partial{t}}(T,t){U}_{2}(\bar{X}^{t,x}(T))\bigg],
$$
with boundary condition $v(T,x)=n {U}_{2}(x),$ and the processes $\bar{X}$ of \eqref{0dyn},
and $ \bar{Z}^{t,x} (s)\triangleq \eta(s)\bar{X}^{t,x} (s)+ l(s) F_3(s, \bar{X}^{t,x} (s)).$
\end{lemma}

We now proceed to the second step. In view of function's $v$ concavity in variable $x,$ and the definition of $(F_1, F_2, F_3)$ (see \eqref{109con}), the equation \eqref{12dE} can be
re-written as
\begin{equation}\label{13ddE}
\!\!\!\!\!\!\!\!\frac{\partial v}{\partial t}(t,x)+\sup_{\zeta ,c,p}\bigg[ \bigg(r+\mu \zeta -c-p+i(t)\bigg)
\frac{\partial v}{\partial x}(t,x)+$$$$\frac{1}{2}\sigma ^{2}\zeta ^{2}%
\frac{\partial ^{2}v}{\partial x^{2}}(t,x)+U_{1}(c)+m{U}_{3}(\eta(t)x+l(t)p)\bigg]  =
\end{equation}
$$
 \mathbb{E}\bigg[
\int_{t}^{T}  \frac{\partial Q}{\partial{t}} (s,t)U_{1}(F_2(s, \bar{X}^{t,x}(s) )  )\,ds+\int_{t}^{T}\frac{\partial q}{\partial{t}} (s,t){U}_{3}(\bar{Z}^{t,x}(s))\,ds+n \frac{\partial Q}{\partial{t}}(T,t){U}_{2}(\bar{X}^{t,x}(T))\bigg].
$$
We notice that

$$
{J}(t,x,\zeta_{\epsilon},c_{\epsilon}, p_{\epsilon})- J(t,x,\bar{\zeta},\bar{c}, \bar{p})=$$$$
 \mathbb{E}\bigg[ \int_{t}^{t+\epsilon} Q(s,t)[U_{1}(c(s))-U_{1}(F_2(s,{X}^{t,x}(s)))]\,ds\bigg]+$$$$
 \!\!\!\!\! \mathbb{E}\bigg[ \int_{t}^{t+\epsilon} q(s,t)[{U}_{3}(Z^{t,x}(s))-{U}_{3}(\bar{Z}^{t,x}(s))]\,ds \bigg]+$$$$
 \!\!\!\!\! \mathbb{E}[ v(t+\epsilon,{X}^{t,x}(t+\epsilon))-v(t+\epsilon,\bar{X}%
^{t,x}(t+\epsilon))]+$$$$
\!\!\!\!\! \mathbb{E}\left[\int_{t+\epsilon}^{T}  [Q(s,t)-Q(s,t-\epsilon)][U_{1}(F_2(s,%
\bar{X}^{t,x}(s)))-U_{1}(F_2(s,{X}^{t,x}(s)))]\,ds\right]+$$$$\!\!\!\!\!
\mathbb{E}\left[\int_{t+\epsilon}^{T} \!\!\![q(s,t)-q(s,t-\epsilon)][{U}_{3}(Z^{t,x}(s))-{U}_{3}(\bar{Z}^{t,x}(s))]\,ds\right]+$$$$\!\!\!\!\!
\mathbb{E}\left[ n  [Q(T,t)-Q(T,t-\epsilon)] [ {U}_{2}(X^{t,x}(T))- {U}_{2}(\bar{X}^{t,x}(T))]\right].
$$

The RHS of this equation has six terms and we will treat each of these
terms separately:

$1.$ In the light of inequality \eqref{189} and the Lebesgue Dominated
Convergence Theorem

\begin{equation*}
{\lim_{\epsilon\downarrow 0}}\frac{ \mathbb{E}\bigg[ \int_{t}^{t+\epsilon} Q(s,t)[U_{1}(c(s))-U_{1}(F_2(s,{X}^{t,x}(s)))]\,ds\bigg]}{\epsilon}$$$$= [U_{1}(c(t))-U_{1}(F_2(t,x))].
\end{equation*}

$2.$ In the light of inequality \eqref{189} and the Lebesgue Dominated
Convergence Theorem

\begin{equation*}
{\lim_{\epsilon\downarrow 0}}\frac{\mathbb{E}\bigg[ \int_{t}^{t+\epsilon}  q(s,t)[{U}_{3}(Z^{t,x}(s))-{U}_{3}(\bar{Z}^{t,x}(s))]\,ds \bigg] }{\epsilon}= m [{U}_{3}(x+p(t)l(t))-{U}_{3}(x+F_3(t,x) l(t))].
\end{equation*}

$3.$ One has

$$
\mathbb{E}[v(t+\epsilon,\bar{X}^{t,x}(t+\epsilon))- v(t+\epsilon,{X}%
^{t,x}(t+\epsilon))]=$$$$\mathbb{E}[ v(t+\epsilon,\bar{X}^{t,x}(t+\epsilon))- v(t,x)]-\mathbb{E}%
[v(t+\epsilon,{X}^{t,x}(t+\epsilon)- v(t,x)].
$$
Moreover
$$
\mathbb{E}[v(t+\epsilon,\bar{X}^{t,x}(t+\epsilon))-  v(t,x)]=\mathbb{E}%
\int_{t}^{t+\epsilon}d[ v(u, \bar{X}^{t,x}(u))].
$$
It\^{o}'s formula yields

\begin{equation*}
{\lim_{\epsilon\downarrow 0}}\frac{\mathbb{E}\int_{t}^{t+\epsilon} d[  v(u, \bar{X}^{t,x}(u))]   }{\epsilon}=$$$$\bigg[
\frac{\partial {v}}{\partial t}(t,x)+\bigg(rx+\mu F_1(t,x) -F_2(t,x)-F_3(t,x) +i(t)  \bigg) \frac{\partial {v}}{\partial x}(t,x)
+ \frac{\sigma^{2} F^{2}_{1}(t,x)} {2}\frac{\partial^{2} v}{\partial x^{2}}(t,x)\bigg].
\end{equation*}

Similarly
$$
{\lim_{\epsilon\downarrow 0}}\frac{ \mathbb{E}[ v(t+\epsilon,{X}^{t,x}(t+\epsilon))-  v(t,x)] }{\epsilon}=
{\lim_{\epsilon\downarrow 0}}\frac{\mathbb{E}\int_{t}^{t+\epsilon} d[   v(u, {X}^{t,x}(u))]  }{\epsilon}=$$$$  \bigg[ \frac{\partial v}{\partial t}(t,x)+ \bigg(r+\mu \zeta(t) -c(t)-p(t)+i(t)\bigg)
\frac{\partial v}{\partial x}(t,x)+\frac{1}{2}\sigma ^{2}\zeta ^{2} (t)%
\frac{\partial ^{2}v}{\partial x^{2}}(t,x) \bigg].
$$

$4.$ In the light of inequality \eqref{189} and the Lebesgue Dominated
Convergence Theorem it follows that

\begin{equation*}
{\lim_{\epsilon\downarrow 0}}\frac{ \mathbb{E}\left[\int_{t+\epsilon}^{T}  [Q(s,t)-Q(s,t-\epsilon)][U_{1}(F_2(s,%
\bar{X}^{t,x}(s)))-U_{1}(F_2(s,{X}^{t,x}(s)))]\,ds\right]}{\epsilon}=0.
\end{equation*}

$5.$ Similarly

\begin{equation*}
{\lim_{\epsilon\downarrow 0}}\frac{\mathbb{E}\left[\int_{t+\epsilon}^{T} [q(s,t)-q(s,t-\epsilon)][{U}_{3}(Z^{t,x}(s))-{U}_{3}(\bar{Z}^{t,x}(s))]\,ds\right]}{\epsilon}=0.
\end{equation*}

$6.$ Finally, by the same token

\begin{equation*}
{\lim_{\epsilon\downarrow 0}}\frac{\mathbb{E}\left[ n  [Q(T,t)-Q(T,t-\epsilon)] [ {U}_{2}(X^{t,x}(T))- {U}_{2}(\bar{X}^{t,x}(T))]\right]}{\epsilon}=0.
\end{equation*}

Therefore

$$
{\lim_{\epsilon\downarrow 0}}\frac{J(t,x,\bar{\zeta},\bar{c}, \bar{p})- {J}(t,x,\zeta_{\epsilon},c_{\epsilon}, p_{\epsilon})}{\epsilon}=
$$$$
\!\!\!\!\!\!\!\!\!\!  \bigg[ \frac{\partial {v}}{\partial t}(t,x)+\bigg(rx+\mu F_1(t,x) -F_2(t,x)-F_3(t,x) +i(t)  \bigg) \frac{\partial {v}}{\partial x}(t,x)+ $$$$
\frac{\sigma^{2} F^{2}_{1}(t,x)} {2}\frac{\partial^{2} v}{\partial x^{2}}(t,x)+U_{1}(F_2(t,x)))+ m {U}_{3}(x+l(t)F_3(t,x)) \bigg]-
$$$$
\!\!\!\!\!\!\!\!\!\! \bigg[ \frac{\partial v}{\partial t}(t,x)+ \bigg(r+\mu \zeta(t) -c(t)-p(t)+i(t))\bigg)
\frac{\partial v}{\partial x}(t,x)+$$$$\frac{1}{2}\sigma ^{2}\zeta ^{2}(t)%
\frac{\partial ^{2}v}{\partial x^{2}}(t,x)+U_{1}(c(t))+m {U}_{3}(x+l(t)p(t))\bigg]\geq 0,
$$

where the last inequality comes from \eqref{12dE} and \eqref{13ddE}.

\begin{flushright}
$\square$
\end{flushright}

\section{CRRA Preferences}
Finding a value function is a complicated problem. We are able to deal with the case of power type
utilities, that is (with some abuse of notations) $U_{1}(x)=U_{2}(x)=U_{3}(x)=U_{\gamma} (x)\triangleq\frac{x^{\gamma}}{\gamma},$ with $\gamma<1.$ In this case we search for the value function $v$ of the form
\begin{equation}
v(t,x)=a(t)U_{\gamma }(x+b(t)),  \label{vo*}
\end{equation}%
where the functions $a(t),$ $b(t)$ are to be found. We consider here
the case $\gamma \neq 0$ (the case of logarithmic utility will be treated
separately). In the light of equations \eqref{109con} one gets
\begin{eqnarray}
F_{1}(t,x) &=&\frac{\mu (x+b(t))}{\sigma ^{2}(1-\gamma )},\,%
\,F_{2}(t,x)=[a(t)]^{\frac{1}{\gamma -1}}(x+b(t)),  \label{0109con} \\
F_{3}(t,x) &=&\frac{1}{l(t)}\left[ ([\frac{a(t)}{m}]^{\frac{1}{\gamma -1}%
}-\eta (t))x+[\frac{a(t)}{m}]^{\frac{1}{\gamma -1}}b(t))\right] .
\label{00109con}
\end{eqnarray}
By \eqref{1o*} the associated wealth process has the following dynamics:
\begin{eqnarray*}
\!\!\!\!\!d\bar{X}(s) &=&\bigg[\left( r+\frac{\eta (s)}{l(s)}\right) \bar{X}%
(s)+\frac{\mu ^{2}}{\sigma ^{2}(1-\gamma )}(\bar{X}(s)+b(s))\\&-&(a(s))^{\frac{1%
}{\gamma -1}}\left( 1+\frac{1}{m^{\frac{1}{\gamma -1}}l(s)}\right) (\bar{X}%
(s)+b(s))\bigg]ds \\
&&+i(s)ds+\frac{\mu (\bar{X}(s)+b(s))}{\sigma (1-\gamma )}dW(s).
\end{eqnarray*}

Let us define the process ${Y}(s)\triangleq \bar{X}(s)+b(s)$ which has the
dynamics
\begin{eqnarray}\label{?!}
d{Y}(s) &=&\bigg[(r+\frac{\eta (s)}{l(s)})Y(s)+\frac{\mu ^{2}}{\sigma
^{2}(1-\gamma )}Y(s)-[a(s)]^{\frac{1}{\gamma -1}}\left( 1+\frac{1}{m^{\frac{1%
}{\gamma -1}}l(s)}\right) Y(s) \\\notag
&&+i(s)+b^{\prime }(s)-(r+\frac{\eta (s)}{l(s)})b(s)\bigg]ds+\frac{\mu Y(s)}{%
\sigma (1-\gamma )}dW(s).
\end{eqnarray}
For considerations that will become clear later on we choose $b(s)$ such
that
\begin{equation*}
i(s)+b^{\prime }(s)-(r+\frac{\eta (s)}{l(s)})b(s)=0,\,\,\,\mbox{and}%
\,\,\,b(T)=0.
\end{equation*}%
By solving this ODE we get
\begin{equation}
b(s)=\int_{s}^{T}i(u)e^{-\int_{u}^{s}\left( r+\frac{\eta (x)}{l(x)}\right)
\,dx}du.  \label{?}
\end{equation}
Solving for the process $Y(s)$ we get

\begin{eqnarray*}
Y(s)&=&Y(t) \exp\bigg(\int_{t}^{s}\bigg(r+\frac{\mu^2}{2\sigma^2 (1-\gamma)^2}+%
\frac{\eta(u)}{l(u)}-(a(u))^{\frac{1}{\gamma-1}} \bigg(1+\frac{1}{m^{\frac{1%
}{\gamma-1}} l(u)}\bigg) \bigg)du\\&+& \frac{\mu(W(s)-W(t))}{\sigma(1-\gamma)}
\bigg)
\end{eqnarray*}

Therefore

\begin{eqnarray*}
\bar{X}^{t,x}(T)&=&x\exp\bigg(\int_{t}^{T}\bigg(r+\frac{\mu^2}{2\sigma^2
(1-\gamma)^2}+\frac{\eta(u)}{l(u)}-(a(u))^{\frac{1}{\gamma-1}} \bigg(1+\frac{%
1}{m^{\frac{1}{\gamma-1}} l(u)}\bigg) \bigg)du\\&+& \frac{\mu(W(T)-W(t))}{%
\sigma(1-\gamma)} \bigg).
\end{eqnarray*}
By plugging $v$ of \eqref{vo*} ( with $(F_{1},F_{2},F_{3})$ of %
\eqref{0109con}, \eqref{00109con}) into \eqref{100ie19(*} and \eqref{1o*}, we
obtain the following integral equation (IE) for $a(t)$
\begin{eqnarray}
a(t)\!\! \!&=&\!\!\! \int_{t}^{T}\!\![Q(s,t)+q(s,t)](a(s))^{\frac{\gamma }{\gamma -1}%
}e^{K(s-t)+\left( \int_{t}^{s}\frac{\gamma \eta (z)}{l(z)}-\gamma (a(z))^{%
\frac{1}{\gamma -1}}\left( 1+\frac{1}{m^{\frac{1}{\gamma -1}}l(z)}\right)
\,dz\right) }\,ds  \label{IE0} \\
&+&nQ(T,t)e^{K(T-t)+\left( \int_{t}^{T}\frac{\gamma \eta (z)}{l(z)}-\gamma
(a(z))^{\frac{1}{\gamma -1}}\left( 1+\frac{1}{m^{\frac{1}{\gamma -1}}l(z)}%
\right) dz\right) },\qquad a(T)=n.  \notag
\end{eqnarray}%
with
\begin{equation}\label{Kk}
K\triangleq \gamma \left( r+\frac{\mu ^{2}}{2(1-\gamma )\sigma ^{2}}\right) .
\end{equation}

Let us summarize this finding:

\begin{lemma}\label{op0}
\label{int} Let $a(t)$ be a solution of the fixed-point problem \eqref{IE0}.
Define $b(t)$ by \eqref{?}. Then $v(t,x)=a(t)U_{\gamma }(x+b(t))$ is a value
function.
\end{lemma}

We turn our attention to the integral equation \eqref{IE0}. Set

\begin{equation}\label{Mm}
M(z)\triangleq 1+\frac{1}{m^{\frac{1}{\gamma -1}}l(z)}
\end{equation}

\begin{assumption}
\label{A1} We require that:
\begin{equation}
\min_{t\in \left[ 0,T\right] }(1-\gamma M(t)+\lambda (t))\geq 0.  \label{&^}
\end{equation}
\end{assumption}

The Assumption \ref{A1} is met if $\gamma \leq 0.$ If $m=1$ and $l(t)=%
\frac{1}{\lambda (t)},$ then  Assumption \ref{A1} is also satisfied (this is
the situation considered by \cite{Pliska}). In the case when $\min_{t\in %
\left[ 0,T\right] }(1-\gamma M(t)+\lambda (t))<0$ it can happen that $a(t)$
reaches 0 which leads to unbounded consumption.

The following Proposition is proved in Appendix C.

\begin{proposition}\label{pio}
\label{existenceODE} If Assumption \ref{A1} is satisfied, then there exists
a unique global $C^{1}$ solution of the integral equation \eqref{IE0}.
\end{proposition}

In other words, for the problem under consideration, there always exists a
value function of the special type $v(t,x)=a(t)U_{\gamma }(x+b(t))$ (note
that there may be others as well). We now proceed to deduce the existence of
policies.

\begin{theorem}\label{existence1}
Let $v$ be the value function of Lemma \ref{op0}. Then, $\{\bar{\zeta}(s),\bar{c}(s),\bar{p}(s)\}_{s\in
\lbrack t,T]}$ given by \eqref{ie} is a policy.
\end{theorem}
The proof follows from Theorem \ref{existence}, Lemma \ref{op0} and Proposition \ref{pio}
as long as we prove that  $\{\bar{\zeta}(s),\bar{c}(s),\bar{p}(s)\}_{s\in
\lbrack t,T]}$ is an admissible strategy and Assumption \ref{A2} is met. The first claim follows
if one establishes \eqref{189}. Taking into account the special form of $v,$ $U_{\gamma},$
and $\bar{X}(s)+b(s),$ (see \eqref{?!}) then Burkholder Davis Gundy inequality yields \eqref{189}. Next,
to show that Assumption \ref{A2} holds true boils down to construct a $C^{1,2}$ solution to some PDEs. Again
by exploiting the special structure, one can construct solutions (for the PDEs of Assumption \ref{A2}) of the form $l(t) U_{\gamma}(x+b(t)),$ with
$l(t)$ being the solution of some ODE.

\subsection{The Case of Logarithmic Utility }

In this special case we can solve the integral equation \eqref{IE0} in closed form.
Indeed with $\gamma=0$ (the case of logarithmic utility) we follow the ansatz
\begin{equation}  \label{vo**}
v(t,x)=a(t)U_{\gamma}(x+b(t))+d(t).
\end{equation}

Then \eqref{IE0} becomes

\begin{equation}\label{I1E0}
a(t)= \int_{t}^{T} [Q(s,t)+q(s,t)]\,ds+n Q(T,t),
\end{equation}
with $b(t)$ given in \eqref{?} and an appropriate choice of function $d(t)$. The equilibrium policy is then given through \eqref{109con} which becomes

\begin{eqnarray}
F_{1}(t,x) &=&\frac{\mu (x+b(t))}{\sigma ^{2}},\,%
\,F_{2}(t,x)=[a(t)]^{-1}(x+b(t)),  \label{10109con} \\
F_{3}(t,x) &=&\frac{1}{l(t)}\left[ ([\frac{a(t)}{m}]^{-1}%
-\eta (t))x+[\frac{a(t)}{m}]^{-1}b(t))\right] .
\label{100109con}
\end{eqnarray}

\begin{remark}
Let us notice that the amount invested in the stock is the same as in the case of the standard Merton problem
with exponential discounting. This somehow surprising result is explained by constant return and volatility
for the stock. We conjecture that in a stochastic volatility model these amounts will be different. The consumption and insurance policies differ from the optimal ones except for the case of exponential discounting. In fact, this is the topic of the next subsection.
\end{remark}

\subsection{The Classical Merton Problem}

The case of exponential discounting, $h(t)=\hat{h}(t)=e^{-\rho t}$ and constant Pareto weight $m(t)=m,$
 deserves special consideration. In that case, the equation
 \eqref{12dE} becomes the classical HJB equation given by dynamic programming
 \begin{equation}  \label{1112dE}
\!\!\!\!\!\!\!\!-(\lambda(t)+\rho){v}(t,x)+\frac{\partial {v}}{\partial t}(t,x)+\bigg(rx+\mu F_1(t,x) -F_2(t,x)-F_3(t,x) +i(t)  \bigg) \frac{\partial {v}}{\partial x}(t,x)+$$$$
\frac{\sigma^{2} F^{2}_{1}(t,x)} {2}\frac{\partial^{2} v}{\partial x^{2}}(t,x)+U_{\gamma}(F_2 (t,x)) +m {U}_{\gamma}(\eta(t)x+l(t)F_3(t,x))=0,
\end{equation}
with the boundary condition $v(T,x)=n {U}_{\gamma}(x),$ and $(F_1, F_2, F_3)$ given through \eqref{109con}. Therefore for the case of exponential discounting the optimal strategy
given by dynamic programming coincides with the policy (given through \eqref{100ie19(*}, \eqref{1o*} \eqref{109con}). This non-linear equation can be linearized by Fenchel-Legendre
transform and therefore it can be shown that it has a unique solution. Moreover, it can be computed by the ansatz \eqref{vo*}. The function $a(t)$ solves an ODE which can be solved explicitly to yield
$$a(t)=\left[  n^{\frac{1}{1-\gamma}}e^{\int_{t}^{T}  \frac{K+\frac{\g\eta(s)}{l(s)}-\rho-\lambda(s)}{1-\gamma} ds} +\int_{t}^{T}
\left( \frac{1+\lambda(u)-\gamma M(u)}{1-\gamma}\right) e^{\int_{t}^{u}  \frac{K+\frac{\g\eta(s)}{l(s)}-\rho-\lambda(s)}{1-\gamma} ds}   du\right]^{1-\gamma},$$

with $K$ given by \eqref{Kk} and $M(z)$ given by \eqref{Mm}.

 \subsection{The Merton Problem with Hyperbolic Discounting}
In this section we assume that the decision-maker gets no income ($i(t)=0$), he/she does not buy life insurance
and there is no possibility of him/her dying before $T.$ Moreover, we assume that discounting is hyperbolic, i.e.,
$h(t)=(1+k_1 t)^{-\frac{k_2}{k_1}},$ with $k_1, k_2$ positive. In \cite{LoPre} it is shown that CRRA type utilities and hyperbolic or
exponential discounting exhibit \textit{the common difference effect}.
Due to this effect, people are more sensitive to a given
time delay if it occurs earlier than later. More precisely, if a person is
indifferent between receiving $x>0$ immediately, and $y>x$ at some later
time $s,$ then he or she will strictly prefer the better outcome if both
outcomes are postponed by some time $t:$
\begin{equation*}
U(x)=h(s)U(y),\qquad \mbox{implies\,\,\,that}\qquad U(x)h(t)<h(t+s)U(y).
\end{equation*}
Furthermore, they assume that the delay needed to compensate for the larger
outcome is a linear function of time, that is
\begin{equation*}
U(x)=h(s)U(y),\qquad \mbox{implies\,\,\,that}\qquad U(x)h(t)=h(kt+s)U(y),
\end{equation*}
for some constant $k.$ They show that the only solution of this functional
equation is $U(x)=\frac{x^{\gamma}}{\gamma}.$
We pay special attention to this case because it explains the consumption puzzle; there is a satiation in the consumption rate
before maturity and exponential discounting can not
capture it (in fact with exponential discounting the optimal consumption rate is either increasing or decreasing at
all times depending on the relationship between the discount rate and the interest rate). Moreover, it says that
optimal strategies and policies are not observationally equivalent. In the following we illustrate this point by a numerical experiment.

We consider one stock
following a geometric Brownian motion with drift $\alpha=0.12,$ volatility $%
\sigma=0.2,$ interest rate $r=0.05,$ and the horizon $T=4.$ This set of parameters is
chosen for illustration. Inspired by \cite{Lai} let the discount
function $h(x)=(1+k_1 x)^{-\frac{k_2}{k_1}}$ be one of the three choices of
hyperbolic discount: case $1.\,\, k_1=5$; case $2.\,\, k_2=10$; case $3.\,\,
k_3=15;$ and $b$ is chosen such that $h(1)=0.3.$ We set  $%
\gamma=-1$ (this choice reflects risk aversion). Let us consider three cases: $n=1, 10, 30.$ We apply the numerical scheme developed in the Numerical Results Section.

\vspace{0.5cm}

As we see from these graphs, the consumption rate policy is
increasing up to a satiation point after which it is decreasing. This phenomena
is observed from the data (people are consuming more and more up to some age
(around 60 years) after which the consumption is decreasing). This may be explained by a drop in income. As the parameter $n$ gets higher (when the agents get more utility from terminal
wealth) the satiation point comes earlier.

 Lemma \ref{13} shows that consumption rate policy is not always monotone.

 \begin{lemma}\label{nonMon}
\label{13}
One can find a hyperbolic discount function such that the consumption rate policy
is neither increasing nor decreasing in time.
\end{lemma}

Appendix D proves this Lemma.\\\\

\subsection{The Stationary Case}

Let us now consider the stationary problem. The coefficients in the model are assumed constant
taking their stationary value, i.e., $n=0,$ $m(t)=m,$ $l(t)=l,$ $\lambda (t)=\lambda,$ $i(t)=i, \eta(t)=\eta,$ and $T=\infty.$ For simplicity
 we assume that
  $$q(s,t)= m \lambda \exp\{-(\lambda+r_1)(s-t)\},\qquad Q(s,t)= \exp\{-(\lambda+r_2)(s-t)\},$$
for some $r_1$ and $r_2$ positive. Before engaging into the formal definition,
let us point the following key fact.  For an admissible time homogeneous (stationary)
policy process $\{{\zeta}(t),{c}(t), p(t)\}_{t\in[0,\infty)}$ and its corresponding wealth process $\{X(t)\}_{t\in[0,\infty)}$
(see \eqref{equ:wealth-one})
the expected utility functional $J(t,x,\zeta,c)$ is time homogeneous, i.e.,
\begin{eqnarray*}
J(t,x,\zeta,c,p)&=&J(0,x,\zeta,c,p) \\\notag &\triangleq&
\mathbb{E}\left[\int_{0}^{\infty}\!\!\!\!\!\!  \exp\{-(\lambda+r_2)s\}U(c(s))\,ds+
\int_{0}^{\infty}\!\!\!\!\!\! m\lambda \exp\{-(\lambda+r_1)s\}U(Z^{0,x}(s))\,ds\right].
\end{eqnarray*}
The intuition is that the clock can be reset so that the expected utility criterion takes zero as
the starting point. We have a similar definition for policies as in the case of finite horizon.

\begin{definition}\label{finite^}

An admissible trading strategy  $\{\bar{\zeta}(s),\bar{c}(s), \bar{p}(s) \}_{s\in \lbrack 0,\infty]}$ is a policy if
there exists a map  $F=(F_{1},F_{2}, F_3):\mathbb{R}\rightarrow \mathbb{R}\times \lbrack 0,\infty )\times \mathbb{R}$  such that for any $x>0$
\begin{equation}
{\lim \inf_{\epsilon \downarrow 0}}\frac{J(0,x,\bar{\zeta},\bar{c}, \bar{p})-J(0,x,\zeta
_{\epsilon },c_{\epsilon }, p_{\epsilon})}{\epsilon }\geq 0,  \label{opt}
\end{equation}%
where

\begin{equation}
\bar{\zeta}(s)={F_{1}(\bar{X}(s))},\quad \bar{c}(s)=%
F_{2}(\bar{X}(s)), \quad  \bar{p}(s)=F_{3}(\bar{X}(s))   \label{^0000eq}
\end{equation}%
and the wealth process $\{\bar{X}(s)\}_{s\in \lbrack 0,\infty]}$ is a solution of
the stochastic differential equation (SDE)
\begin{equation}
d\bar{X}(s)=[r\bar{X}(s)+\mu F_{1}(\bar{X}(s))-F_{2}(\bar{X}(s))-F_{3}(\bar{X}(s))+i(s)]ds+\sigma F_{1}(\bar{X}%
(s))dW(s).  \label{^0dyn}
\end{equation}

Moreover, the process $\{{\zeta}_{\epsilon}(s),{c}_{\epsilon}(s),  p_{\epsilon}(s)\}_{s\in[0,\infty]}$ is
another investment-consumption strategy defined by
\begin{equation}  \label{1e}
\zeta_{\epsilon}(s)=%
\begin{cases}
\bar{\zeta}(s),\quad s\in[0,\infty]\backslash E_{\epsilon} \\
\zeta(s), \quad s\in E_{\epsilon},%
\end{cases}%
\end{equation}

\begin{equation}  \label{2e}
c_{\epsilon}(s)=%
\begin{cases}
\bar{c}(s),\quad s\in[0,\infty]\backslash E_{\epsilon} \\
c(s), \quad s\in E_{\epsilon},%
\end{cases}%
\end{equation}
\begin{equation}  \label{3e}
p_{\epsilon}(s)=%
\begin{cases}
\bar{p}(s),\quad s\in[0,\infty]\backslash E_{\epsilon} \\
p(s), \quad s\in E_{\epsilon},%
\end{cases}%
\end{equation}
with $E_{\epsilon}=[0,\epsilon],$ and $\{{\zeta}(s),{c}(s),{p}(s)\}_{s\in
E_{\epsilon} }$ is any strategy for which $\{{\zeta}_{\epsilon}(s),{c}%
_{\epsilon} (s), {p}_{\epsilon} (s)\}_{s\in[0,\infty]}$ is an admissible policy.
\end{definition}

Similarly we define the value function.
\begin{definition}\label{de1^}
 A function $v:\mathbb{R}\rightarrow  \mathbb{R}$ is a value function if
it solves the following system of equations

\begin{equation}
v(x)=J(0,x,\bar{\zeta},\bar{c}, \bar{p})
\label{00ie19(*}
\end{equation}%

$$\bar{\zeta}(s)={F_{1}(\bar{X}(s))},\quad \bar{c}(s)=%
F_{2}(\bar{X}(s)), \quad  \bar{p}(s)=F_{3}(\bar{X}(s)), $$

\begin{equation}\label{o*}
d\bar{X}(s)=[r\bar{X}(s)+\mu F_{1}(\bar{X}(s))-F_{2}(\bar{X}%
(s))- F_{3}(\bar{X}(s))+i(s)]ds+\sigma F_{1}(\bar{X}(s))dW(s),\end{equation}%

\begin{equation}\label{))+}
F_{1}(x)=-\frac{\mu v'(x)}{\sigma ^{2}
v''(x)},\,\,F_{2}(x)=\left(v'(x))\right)^{\frac{1}{\gamma-1}} ,\,\,F_{3}(x)=\frac{1}{l} \left[ \left( \frac{1}{m}v'(x)\right)^{\frac{1}{\gamma-1}}-\eta x\right].
\end{equation}
\end{definition}

Let us look for the value function of the form
\begin{equation}  \label{vo!}
v(x)=aU_{\gamma}(x+b),
\end{equation}

for some constants $a$ and $b.$  By solving $(3.4)$, we get $b=\frac{i}{r+\frac{\eta}{l}}$.
Let $\beta=1+\frac{m^{\frac{1}{1-\g}}}{l}$and $K$ be given by \eqref{Kk}. The constant $a$ should solve the following equation
 \begin{equation}\label{eq_a}
a=\frac{a^{\frac{\g}{\g-1}}} {\lam+r_1-K-\frac{\g\eta}{l}+\g\beta a^{\frac{1}{\g-1}}} + m \lam \frac{  (\frac{a}{m})^{\frac{\g}{\g-1}}  } {\lam+r_2-K-\frac{\g\eta}{l}+\g\beta a^{\frac{1}{\g-1}}  }
\end{equation}  with the transversality conditions
\begin{equation}\label{TC}
\lam+r_j-K-\frac{\g\eta}{l}+\g\beta a^{\frac{1}{\g-1}}>0\qquad j=1,2.\end{equation}

\begin{lemma}\label{QuadraticSol}
There is a unique solution of \eqref{eq_a} and \eqref{TC}.
\end{lemma}

 Appendix E proves this Lemma.\\\\

 We are ready to state the main result of this section.

\begin{theorem}\label{existence1}
Let $v$ be defined by \eqref{vo!} with $a$ the solution of \eqref{eq_a}. The function $(F_1, F_2, F_3)$ of \eqref{))+} defines a policy
through \eqref{^0dyn} and \eqref{^0000eq}.
\end{theorem}

Proof: The proof for the most part parallels Theorem \ref{existence}. The only part which requires more analysis is
showing that

\begin{equation*}
{\lim_{\epsilon\downarrow 0}}\mathbb{E}\left[\int_{\epsilon}^{\infty}  \exp\{-(\lambda+r_i)s\}[U_{\gamma}(F_2(%
\bar{X}^{0,x}(s)))-U_{\gamma}(F_2({X}^{0,x}(s)))]\,ds\right]=0,\quad i=1,2.
\end{equation*}
which is equivalent to

\begin{equation*}
{\lim_{\epsilon\downarrow 0}}\mathbb{E}\left[\int_{\epsilon}^{\infty}   \exp\{-(\lambda+r_i)s\}[(\bar{X}^{0,x}(s)+b)^{\gamma}-({X}^{0,x}(s)+b)^{\gamma}]\,ds\right]=0, ,\quad i=1,2.
\end{equation*}
The result follows from Lebesque Dominated Convergence Theorem if we prove that

\begin{equation}\label{000}
\mathbb{E}\left[\int_{0}^{\infty}   \exp\{-(\lambda+r_i)s\}[(\bar{X}^{0,x}(s)+b)^{\gamma}+({X}^{0,x}(s)+b)^{\gamma}]\,ds\right]<\infty,\quad i=1,2.
\end{equation}
Notice that from the transversality conditions \eqref{TC} one gets

\begin{equation*}
\mathbb{E}\left[\int_{0}^{\infty}   \exp\{-(\lambda+r_i)s\}(\bar{X}^{0,x}(s)+b)^{\gamma}\right]<\infty,\quad i=1,2.
\end{equation*}

Moreover, if $s\in[\epsilon, \infty]$ then

$$\left( \frac{{X}^{0,x}(s)+b}{\bar{X}^{0,x}(s)+b}\right)^{\gamma} =\left( \frac{{X}^{0,x}(\epsilon)+b}{\bar{X}^{0,x}(\epsilon)+b}\right)^{\gamma} \triangleq R(\epsilon),$$
and $R(\epsilon)$ and $\bar{X}^{0,x}(s)$ are independent. Thus, Holder inequality and a standard argument yields
\begin{equation*}
\mathbb{E}\left[\int_{0}^{\infty}   \exp\{-(\lambda+r_i)s\}({X}^{0,x}(s)+b)^{\gamma}\right]<\infty,\quad i=1,2,
\end{equation*}
so \eqref{000} holds true.

\begin{flushright}
$\square$
\end{flushright}

\section{Numerical Results}

We provide a numerical scheme to approximate the integral equation \eqref{IE0}. For simplicity we assume that $\eta(t)=1.$
 In a first step let us discretize the interval $\left[0, T \right]$ by introducing the points $t_n = T + nh$ , where
 $\e=-\frac{T}{N}$; recall that with $K$ given by \eqref{Kk} and $M(\cdot)$ of \eqref{Mm}, the equation \eqref{IE0} becomes
 in a differential form

 \begin{eqnarray}\label{lll}
a'(t) &=&(\g M(t)-\lambda(t)-1)(a(t))^{\frac{\g}{\g-1}}+\left(\lambda(t)-\frac{h'(T-t)}{h(T-t)}-K-\frac{\g}{l(t)}\right)a(t)\\\notag
&+& \int_{t}^{T}L(s,t)(a(s))^{\frac{\g}{\g-1}}\left(\frac{A(s)}{A(t)}\right)ds,
\end{eqnarray}
where $$A(s)\triangleq\exp(\int_s^T \g (a(z))^{\frac{1}{\g-1}}M(z)dz)$$
 and
  $$L(s,t)\triangleq\left[\left(\frac{h'(T-t)}{h(T-t)}-\frac{h'(s-t}{h(s-t)}\right)Q(s,t)+\left(\frac{h'(T-t)}{h(T-t)}-\frac{\bar{h}'(s-t)}{\bar{h}(s-t)}\right)q(s,t)\right]e^{\int_{t}^{s}K+\frac{\g}{l(u)}du}.$$
From the definition of $A(s),$ it follows that

  $$A'(s)=-\g a(s))^{\frac{1}{\g-1}}M(s) A(s).$$

Our approximation scheme is done in three steps. In a first step, we construct the sequence $a_{n}^1$ and $A_{n}^1$
recursively by

$$ a_{n+1}^1\triangleq a_n^1 + \e a'(t_n),\qquad A_{n+1}^1\triangleq A_n^1 + \e A'(t_n). $$

 \begin{lemma}\label{Le1}
 If $a_n^1$ and $A_n^1$, $n = 0\cdots N$ are defined by $a_0^1=1,$ $A^1_0=1$ and
 \\\\

   $\left\{\begin{array}{lll}
a_{n+1}^1&=&a_n^1+\e\big((\g M(t_n)-\lambda(t_n)-1)(a_n^1)^{\frac{\g}{\g-1}}+\left(\lambda(t_n)-\frac{h'(T-t_n)}{h(T-t_n)}-K-\frac{\g}{l(t_n)}\right)a_n^1\\
&+& \int_{t_n}^{T}L(s,t_n)(a(s))^{\frac{\g}{\g-1}} \left(\frac{A(s)}{A(t_n)}\right)ds \\
A_{n+1}^1&=&A_n^1-\g\e (a(t_n))^{\frac{1}{\g-1}}M(t_n) A_n^1\\
\end{array}\right. $

Then there exists a constant $C$ such that

$$|a_n^1-a(t_n)|\leq C|\e| \,\,\mbox{and}\,\,|A_n^1-A(t_n)| \leq C|\e|,\quad \forall n\in   0,1,\cdots, N.$$

 \end{lemma}

 Appendix F proves this Lemma.\\\\

In a second step we discretize the integral $ \int_{t_n}^{T}L(s,t_n)(a(s))^{\frac{\g}{\g-1}}\left(\frac{A(s)}{A(t_n)}\right)ds.$
 This will lead to the following Lemma.

 \begin{lemma}\label{Le2}
  If $a_n^2$ and $A_n^2$, $n = 0\cdots N$ are defined by $a_0^2=1, A_0^2=1$ and\\\\

   $\left\{\begin{array}{lll}
a_{n+1}^2&=&a_n^2+\e(\g M(t_n)-\lambda(t_n)-1)(a_n^2)^{\frac{\g}{\g-1}}+\e\left(\lambda(t_n)-\frac{h'(T-t_n)}{h(T-t_n)}-K-\frac{\g}{l(t_n)}\right)a_n^2\\
&-&\e^2 \sum_{j=0}^{n-1}L(t_j,t_n)(a(t_j))^{\frac{\g}{\g-1}}\left(\frac{A(t_j)}{A(t_n)}\right)\\
A_{n+1}^2&=&A_n^2-\g\e (a_n^2)^{\frac{1}{\g-1}}M(t_n) A_n^2
\end{array}\right. $

Then there exists a constant $C$ such that

$$|a_n^2-a_n^1)|\leq C|\e| \,\,\mbox{and}\,\,|A_n^2-A_n^1| \leq C|\e|,\quad \forall n\in   0,1,\cdots, N.$$

 \end{lemma}

  Appendix G proves this Lemma.\\\\

In a third step we introduce an explicit scheme.

\begin{lemma}\label{Le3}
  If $a_n^3$ and $A_n^3$, $n = 0\cdots N$ are defined by $a_0^3=1, A_0^3=1$ and\\\\
   $\left\{\begin{array}{lll}
a_{n+1}^3&=&a_n^3+\e(\g M(t_n)-\lambda(t_n)-1)(a_n^3)^{\frac{\g}{\g-1}}
+\e\left(\lambda(t_n)-\frac{h'(T-t_n)}{h(T-t_n)}-K-\frac{\g}{l(t_n)}\right)a_n^3\\
&-&\e^2 \sum_{j=0}^{n-1}L(t_j,t_n)(a_j^3)^{\frac{\g}{\g-1}}\left(\frac{A_j^3}{A_n^3}\right)\\
A_{n+1}^3&=&A_n^3-\g\e (a_n^3)^{\frac{1}{\g-1}}M(t_n) A_n^3\\
\end{array}\right.$

Then there exists a constant $C$ such that

$$|a_n^3-a_n^2)|\leq C|\e| \,\,\mbox{and}\,\,|A_n^3-A_n^2| \leq C|\e|,\quad \forall n\in   0,1,\cdots, N.$$

\end{lemma}

  Appendix H proves this Lemma.\\\\

By using the preceding lemmas and Lipschitz continuity of function $a(t)$ we summarize the results of this section by the following
Theorem.

\begin{theorem}\label{th}
Let $a_N(t)$ be the function obtained by the linear interpolation of the points\\ $(t_n=T-\frac{nT}{N}, a_n^3).$  Then
 $$|a_N(t)-a(t)|\leq C|\e|, \qquad \forall  t \in \left[0,T\right],$$
 for some positive constant $C$ independent of $N.$
\end{theorem}

In the following we perform a numerical experiment. Let
$T=4, r=0.05, \mu=0.07,  \sigma=0.2, p=-1, N=1000, \rho=0.8, \lambda(t)=\frac{1}{200}+\frac{9}{8000}t,  l(t)=\frac{1}{\lambda(t)}, \eta(t)=1,$
The discount function is exponential $h(t)=\hat{h}(t)=\exp(-\rho t)$ with $\rho=0.8$.
The Pareto weight is $m(t)=\log(\frac{T+\epsilon-t}{\epsilon})$ with $\epsilon=10^{-15}.$
We choose this set of parameters for illustration. As people get older, perhaps they weigh more their heirs utility; so a time decreasing  aggregation rate seems the natural choice. It is for this reason that we consider a decreasing function $m.$ We plot
the maps $F_2$ and $F_3$ which lead to the policies. Furthermore, we plot the difference in $F_3$ when $m$ is variable as opposed to $m$ is constant. The results show that higher utility weight $m$ leads to higher amount spent on life insurance.

\section{Conclusion and future research}

We have studied the portfolio management problem in which an agent invests in a risky asset,
consumes, and buys life insurance in order to maximize utility of his/her and heirs. Different discount rates for the agent and heirs lead to time inconsistency. Moreover a time
changing aggregation weight will lead to time inconsistency as well. The way we deal with this predicament is
by looking for subgame perfect Nash equilibrium strategies that we call policies.  We find them in special cases.
Our model is rich enough to capture different aspects in portfolio theory. We perform numerical experiments in order
to explain the effect of discounting in one hand and the effect of aggregation on a different hand. Hyperbolic discounting
is emphasized in a Merton type problem (for simplicity we shut off some parameters). The surprising result is that the policies
and optimal strategies are not always observationally equivalent. For example, in certain cases, watching
one's consumption rate we can infer (from its time monotonicity) wether or not is a optimal or a equilibrium one. Indeed, a non monotone consumption rate can not be optimal for the case of exponential discounting.
The consumption rate policy exhibit a satiation point (it has a  hump shaped
behaviour) and this may explain the consumption puzzle. A time varying aggregation rate is benchmarked to a constant one. Our numerical experiment
comes to support the intuition that the more the manager cares for his/her heirs, the more he/she will pay on life insurance.

We have introduced a system of equations: the integral equation (\ref{100ie19(*}) together with the SDE \eqref{1o*} and PDE \eqref{109con}. Their validity has been established in the case when the utility
function and the bequest function are of (CRRA) type, but we think that it
extends to any concave utilities, as in the deterministic case. This paper
can be seen as a first step in the general direction of extending stochastic
control away from the optimization paradigm towards time-consistent
strategies. The mathematical difficulties are considerable: we have no
general existence nor uniqueness theory for the equations (\ref{100ie19(*}) or (%
\ref{12dE}), which replace the classical HJB equation of optimal control. In the present paper we sidestep the difficulty by using an Ansatz, but we hope that future work, by ourselves or
others, will solve these problems.

We conclude by pointing out that most portfolios are not managed by an individual, but by a group, either a professional management team, or the investor himself and his family. As we mentioned in the introduction, we should then introduce one utility function, one psychological discount rate, and one Pareto weight for each member of the group. Consider for instance a group with two members (husband and wife). Member $i$ (for $i=1,2$) has utility $%
u_{i}\left( c_{1},c_{2}\right) $, where $c_{1}$ is the consumption of the
husband and $c_{2}$ is the \ consumption of the wife, and discount factor $%
h_{i}\left( t\right) $, so the utility derived at time $t$ by member $i$
from the couple consuming $\left( c_{1},c_{2}\right) $ at time $s>t$ is $%
h_{i}\left( s-t\right) u_{i}\left( c_{1},c_{2}\right) $. The utilities of
the husband and the wife have to be aggregated by Pareto weights. If for
instance we assume that, as is the case in most couples, one member
specializes in long-term decisions and the other in short-term ones, we
find that there should exist some decreasing function $\mu \left( t\right) $%
, with $0\leq \mu \left( t\right) \leq 1$, such that the behaviour of the
couple between $t$ and $T$ is adequately described by maximising the
intertemporal criterion:
\begin{equation*}
\int_{t}^{T}\left[ \mu \left( s-t\right) h_{1}\left( s-t\right) u_{1}\left(
c_{1}\left( s\right) ,c_{2}\left( s\right) \right) +\left( 1-\mu \left(
s-t\right) \right) h_{2}\left( s-t\right) u_{2}\left( c_{1}\left( s\right)
,c_{2}\left( s\right) \right) \right] ds
\end{equation*}%
plus some terminal criterion (legacy) at time $T$. Even in the case when
both members of the group have a constant psychological discount rate, so
that $h_{i}\left( t\right) =\exp \left( -r_{i}t\right) $, and even if $%
r_{1}=r_{2}$, the group will exhibit time-inconsistency.

Our model covers the particular case when $u_{1}=u_{2}=U.$ This group is time-inconsistent: their expected intertemporal utility is:
\begin{eqnarray*}
J(t,x,\zeta ,c_1, c_2)\!\!\!\!\!&=\!\!\!\!\!& \mathbb{E}\bigg[
\int_{t}^{T} \!\!\! h_1(s-t) U(c_1(s))\,ds\!+\!\!\int_{t}^{T} \!\!\! m(s-t)\,\, h_2(s-t){U}(c_2(s))\,ds\!\\&+&\!h_1(T-t) U(X^{\zeta, c_1, c_2 }(T)) \bigg],
\end{eqnarray*}
which falls within our model. Assuming that the function
$m(\cdot)$ is decreasing with $m(0)\simeq \infty,$ and $m(T)\simeq 0$ would capture the situation when one member plans for short time and the other plans for long time.

The more general case when $u_{1} \neq u_{2}$ together with other macroeconomic problems with heterogeneous agents will be the subject of future research.

\section{Appendix}

{\bf {Appendix A}}: Proof of Lemma \ref{L1}: We first establish that

\begin{equation}\label{q1}
 \mathbb{E}\bigg[
\int_{t}^{T\wedge\tau}h(s-t)U_{1}(c(s))\,ds\bigg]= \mathbb{E}\bigg[
\int_{t}^{T} Q(s,t) U_{1}(c(s))\,ds\bigg]
\end{equation}

In the light of equations \eqref{*0}, \eqref{*1} and the random variable $\tau$ being independent of Brownian motion $W$
it follows that

$$ \mathbb{E}\bigg[
\int_{t}^{T\wedge\tau}h(s-t)U_{1}(c(s))\,ds\bigg]= \mathbb{E}\bigg[
\exp\{-\int_{t}^{T}\lambda(z)\,dz\}\int_{t}^{T}h(s-t)U_{1}(c(s))\,ds+$$$$
\int_{t}^{T}\lambda(u)\exp\{-\int_{t}^{u}\lambda(z)\,dz\} \int_{t}^{u}h(s-t)U_{1}(c(s))\,ds du\bigg]. $$

Moreover

$$\mathbb{E}\bigg[ \int_{t}^{T}\lambda(u)\exp\{-\int_{t}^{u}\lambda(z)\,dz\} \int_{t}^{u}h(s-t)U_{1}(c(s))\,ds du\bigg]= $$
$$-\mathbb{E}\bigg[ \frac{\partial}{\partial u} \left (\exp\{-\int_{t}^{u}\lambda(z)\,dz\} \right) \int_{t}^{u}h(s-t)U_{1}(c(s))\,ds du\bigg],$$

and integration by parts will lead to \eqref{q1}. It is easy to see that

\begin{equation}\label{q2}
\mathbb{E}\bigg[ \bar{h}(\tau-t) {U}_{3}(Z^{t,x}(\tau)) 1_{\{\tau\leq T|\tau>t\}}=\mathbb{E}\bigg[\int_{t}^{T} q(s,t){U}_{3}(Z^{t,x}(s))\,ds\bigg].
\end{equation}

Finally let us prove that

\begin{equation}\label{q3}
 \mathbb{E}\bigg[ nh(\tau-t){U}_{2}(X^{t,x}(T)) 1_{\{\tau>T|\tau>t\}}\bigg]=
 \mathbb{E}\bigg[ n Q(T,t) {U}_{2}(X^{t,x}(T))\bigg].
\end{equation}

This follows from \eqref{*1}. In the light of \eqref{q1} \eqref{q2} and \eqref{q3}, it follows that \eqref{w9} holds true.

\begin{flushright}
$\square$
\end{flushright}

{\bf {Appendix B}}: Proof of Lemma \ref{pDe}:
 Let the functions $(t,x)\rightarrow f^i(t,\cdot,x)$ satisfy the following PDEs:
\begin{equation}\label{PDE1}
\frac{\partial f^i}{\partial t}+(r x+\mu F_{1}(t,x)-F_{2}(t,x)-F_{3}(t,x)+i(t))\frac{\partial f^i}{\partial x}+\frac{\sigma^{2} F_{1}^{2}(t,x)}{2} \frac{\partial^2 f^i}{\partial x^2}=0,
\quad i=1,2,3.
\end{equation}
with the boundary conditions
$$ f^1(t,s,x)=U_{1}(F_{2}(s,x)),\quad  f^2(t,s,x)=U_{3}(\eta (s)x+l(s)F_{2}(s,x)),\quad  f^3(t,T,x)=U_{2}(x).$$
In the light of Assumption \ref{A2}, these PDEs have $C^{1,2}$ solutions. According to the Feyman-Kac's formula
$$f^1 (t,s,x)=\mathbb{E}[U_{1}(F_{2}(s,\bar{X}^{t,x}(s)))],\,f^2 (t,s,x)=\mathbb{E}[U_{3}(\bar{Z}^{t,x}(s))],\,
f^3 (t,T,x)=\mathbb{E}[U_{2}(\bar{X}^{t,x}(T))].$$
Therefore
\begin{equation}\label{1ie1}
v(t,x)=\int_{t}^{T}(Q(s,t) f^1(t,s,x)+ q(s,t) f^2(t,s,x)) \,ds+ n Q(T,t) f^3(t,T,x)
\end{equation}
By differentiating under the integral sign in \eqref{1ie1} we obtain

 \begin{eqnarray}\label{*1ie1}
\frac{\partial v}{\partial t}(t,x)&=&\int_{t}^{T}(Q(s,t) \frac{\partial f^1}{\partial t}(t,s,x)+ q(s,t) \frac{\partial f^2}{\partial t}(t,s,x)) \,ds+ n Q(T,t) \frac{\partial f^3}{\partial t}(t,T,x)
\\\notag&+&(Q(t,t) f^1(t,t,x)+ q(t,t) f^2(t,t,x)) \\\notag &+&
\int_{t}^{T}(\frac{\partial Q}{\partial{t}} (s,t) f^1(t,s,x)+ \frac{\partial q}{\partial{t}} (s,t) f^2(t,s,x)) \,ds+ n \frac{\partial Q}{\partial{t}} (T,t)f^3(t,T,x) ,
\end{eqnarray}

 \begin{equation}\label{*31ie1}
\frac{\partial v}{\partial x}(t,x)=\int_{t}^{T}(Q(s,t) \frac{\partial f^1}{\partial x}(t,s,x)+ q(s,t) \frac{\partial f^2}{\partial x}(t,s,x)) \,ds+ n Q(T,t) \frac{\partial f^3}{\partial x}(t,T,x),
\end{equation}

and

 \begin{equation}\label{*14ie1}
\frac{\partial^2 v}{\partial x^2}(t,x)=\int_{t}^{T}(Q(s,t) \frac{\partial^2 f^1}{\partial x^2}(t,s,x)+ q(s,t) \frac{\partial^2 f^2}{\partial x^2}(t,s,x)) \,ds+ n Q(T,t) \frac{\partial^2 f^3}{\partial x^2}(t,T,x)
\end{equation}

By combining \eqref{*1ie1},  \eqref{*31ie1},  \eqref{*14ie1} and  \eqref{PDE1}, we obtain \eqref{12dE}.

\begin{flushright}
$\square$
\end{flushright}

{\bf {Appendix C}}: Proof of Proposition \ref{existenceODE}. We proceed in a couple of steps.
In a first step, we obtain lower and upper bounds for $a(t).$ The second
step shows that for the case of discount functions being a linear combination of exponentials
the equation becomes an  ODE system for which we have local existence. The local solution together with the bounds obtained
in the first step will lead to global existence. In a last step, we approximate the  discount functions by linear combination of exponentials
and the solution is constructed by the limit of the ODE systems solutions. In the following, we will go over step one only since the second
step follow as in  \cite{EkePir} and the last step follows from a density argument. For simplicity we assume that $\eta(t)=1.$ Now, let us define
$$\bar{a}(t)=a(t) e^{\left[K(t-T)+\int_{t}^{T}\left(-\frac{\g}{l(z)}+\g a(z)^{\frac{1}{\g-1}}M(z)\right)dz\right]}.$$

It follows that

\begin{eqnarray*}
\bar{a}(t)&=&\int_{t}^{T}\left[Q(s,t)+q(s,t)\right][a(s)]^{\frac{\g}{\g-1}}e^{\left[K(s-T)-\int_{s}^{T}\left(\frac{\g}{l(z)}+\g (a(z))^{\frac{1}{\g-1}}M(z)dz\right)\right]}ds+nQ(T,t)
\end{eqnarray*}

Consequently

 $$\bar{a}'(t)=-(Q(t,t)+q(t,t))(a(t))^{\frac{\g}{\g-1}}e^{\left[K(t-T)-\int_{t}^{T}\left(\frac{\g}{l(z)}+\g (a(z))^{\frac{1}{\g-1}}M(z)\right)dz\right]}ds+
n \frac{\partial}{\partial{t}}Q(T,t)$$$$+\int_{t}^{T}\frac{\partial}{\partial{t}}\left[Q(s,t)+ q(s,t)\right](a(s))^{\frac{\g}{\g-1}}e^{\left[K(s-T)+\int_{s}^{T}\left(\frac{\g}{l(z)}+\g (a(z))^{\frac{1}{\g-1}}M(z)\right)dz\right]}ds
$$

From direct calculations
\begin{eqnarray*}
\frac{\partial}{\partial{t}}Q(T,t)&=&(\lambda(t)h(T-t)-h'(T-t))exp\left[-\int_{t}^{T}\lambda(z)dz\right]\\
               &=& (\lambda(t)-\frac{h'}{h}(T-t))Q(T,t)
\end{eqnarray*}

and
\begin{eqnarray*}
\frac{\partial}{\partial{t}}\left[Q(s,t)+q(s,t)\right]&=&(\lambda(t)h(s-t)-h'(s-t))exp\left[-\int_{t}^{s}\lambda(z)dz\right]\\
&+&\lambda(s)(\lambda(t)\bar{h}(s-t)-\bar{h}'(s-t))exp\left[-\int_{t}^{s}\lambda(z)dz\right].
\end{eqnarray*}
Therefore

  $$ \frac{\partial}{\partial{t}}\left[Q(s,t)+ q(s,t)\right] =\left(\lambda(t)-\frac{h'(s-t)}{h(s-t)}(s-t)\right)Q(s,t)+\left(\lambda(t)-\frac{\bar{h}'(s-t)}{\bar{h}(s-t)}\right) q(s,t)
$$

Thus

\begin{eqnarray*}
\bar{a}'(t)&=&-(\lambda(t)+1)(a(t))^{\frac{\g}{\g-1}}e^{\left[K(t-T)-\int_{t}^{T}\left(\frac{\g}{l(z)}+\g [a(z)]^{\frac{1}{\g-1}}M(z)\right)dz\right]}ds\\
&+&n\left(\lambda(t)-\frac{h'(T-t)}{h(T-t)}\right)Q(T,t)\\
&+& \int_{t}^{T}\left[\left(\lambda(t)-\frac{h'(s-t)}{h(s-t)}\right)Q(s,t)+\left(\lambda(t)-\frac{\bar{h}'(s-t)}{\bar{h}(s-t)}\right)q(s,t) \right]
\\&\times&(a(s))^{\frac{\g}{\g-1}}e^{\left[K(s-T)+\int_{s}^{T}\left(\frac{\g}{l(z)}+\g (a(z))^{\frac{1}{\g-1}}M(z)\right)dz\right]}ds
\end{eqnarray*}

Hence

 \begin{eqnarray*}
\big((a'(t)&+&(K+\frac{\g}{l(t)}-\g (a(t))^{\frac{1}{\g-1}}M(t))a(t)\big)e^{\left[-\int_{t}^{T}\left(K+\frac{\g}{l(z)}\right)dz+\int_{t}^{T}\g (a(z))^{\frac{1}{\g-1}}M(z)dz\right]}\\
&=&-(\lambda(t)+1)(a(t))^{\frac{\g}{\g-1}}e^{\left[-\int_{t}^{T}\left(K+\frac{\g}{l(z)}\right)dz+\int_{t}^{T}\g (a(z))^{\frac{1}{\g-1}}M(z)dz\right]}\\
&+&n\left(\lambda(t)-\frac{h'(T-t)}{h(T-t)})\right)Q(T,t)\\
&+& \int_{t}^{T}\left[\left(\lambda(t)-\frac{h'(s-t)}{h(s-t)}\right)Q(s,t)+\left(\lambda(t)-\frac{\bar{h}'(s-t)}{\bar{h}(s-t)}\right)q(s,t) \right]
\\&\times&(a(s))^{\frac{\g}{\g-1}}e^{\left[-\int_{s}^{T}\left(K+\frac{\g}{l(z)}\right)dz+\int_{s}^{T}\g (a(z))^{\frac{1}{\g-1}}M(z)dz\right]}ds
\end{eqnarray*}

Consequently

$$
\left[(a'(t)+(K+\frac{\g}{l(t)}-\g (a(t))^{\frac{1}{\g-1}}M(t))a(t)\right] =$$$$-(\lambda(t)+1)(a(t))^{\frac{\g}{\g-1}}
+n\left(\lambda(t)-\frac{h'(T-t)}{h(T-t)}\right)Q(T,t)e^{\left[\int_{t}^{T}K+\frac{\g}{l(z)}dz -\int_{t}^{T}\g (a(z))^{\frac{1}{\g-1}}M(z)dz\right]}ds+$$$$
\int_{t}^{T}\left[\left(\lambda(t)-\frac{h'(s-t)}{h(s-t)}\right)Q(s,t)+\left(\lambda(t)-\frac{\bar{h}'(s-t)}{\bar{h}(s-t)}\right)q(s,t) \right]
(a(s))^{\frac{\g}{\g-1}}e^{\left[\int_{t}^{s}\left(K+\frac{\g}{l(z)}\right)dz+\int_{s}^{t}\g (a(z))^{\frac{1}{\g-1}}M(z)dz\right]}ds$$

From the definition of $a(t)$ we get that

\begin{equation}\label{211}
\left[(a'(t)+(K+\frac{\g}{l(t)}-\g (a(t))^{\frac{1}{\g-1}}M(t))a(t)\right] =-(\lambda(t)+1)(a(t))^{\frac{\g}{\g-1}}
+\lambda(t)a(t)-\frac{h'(T-t)}{h(T-t)}a(t)
\end{equation}
$$+ \int_{t}^{T}\bigg(\frac{h'(T-t}{h(T-t})\big(Q(s,t)+q(s,t)\big)-\frac{h'(s-t)}{h(s-t)}(Q(s,t)-$$$$
\frac{\bar{h}'(s-t)}{\bar{h}(s-t)})q(s,t) \bigg)
(a(s))^{\frac{\g}{\g-1}}e^{\left[K(s-t)+\int_{s}^{t}\g [a(z)]^{\frac{1}{\g-1}}M(z)dz\right]}ds
$$

Since $-\rho\leq \frac{h'}{h}(z)\leq \rho$ and $-\rho\leq \frac{\bar{h}'}{\bar{h}}(z)\leq \rho$ and $-\rho' \leq \frac{\g}{l(z)}\leq \rho' $  for $0\leq z\leq T$
the equation \eqref{211} leads to

\begin{equation}\label{)}
a'(t)\leq -(1+\lambda(t)-\g M(t))(a(t))^{\frac{\g}{\g-1}}+(\lambda(t)+3\rho-K+\rho')a(t)
\end{equation}

and

\begin{equation}\label{))}
a'(t)\geq -(1+\lambda(t))(a(t))^{\frac{\g}{\g-1}}-(K+\lambda(t)+3\rho-\rho')a(t)
\end{equation}

Let us denote
$$C_1\triangleq\min_{t\in\left[0,T\right]} (1-\g M(t)+\lambda(t)),\qquad C_0\triangleq\max_{t\in\left[0,T\right]} (\lambda(t)+3\rho-K+\rho')$$

$$D_0\triangleq\max_{t\in\left[0,T\right]} (K+\lambda(t)+3\rho),\qquad D_1\triangleq\max_{t\in\left[0,T\right]}(1+\lambda(t))$$
so that equations \eqref{)} and \eqref{))} will become

\begin{equation}\label{pp1}
a'(t)\leq -C_1(a(t))^{\frac{\g}{\g-1}}+C_0 a(t)
\end{equation}

 and
\begin{equation}\label{pp}
a'(t)\geq -D_1 (a(t))^{\frac{\g}{\g-1}} -D_0 a(t)
\end{equation}

Now, $C_0,D_1,D_0>0$ and by assumption \ref{A1}, it follows that $C_1>0$. Consequently, we get lower and upper bounds on $a(t)$ by integrating \eqref{pp1} and \eqref{pp}  as in  \cite{EkePir}.

\begin{flushright}
$\square$
\end{flushright}

{\bf {Appendix D}}: Proof of Lemma \ref{nonMon}: We take $n=2$ and the coefficients
as in the numerical experiment.  In the light of the differential equation

\begin{eqnarray}  \label{&1234}
a^{\prime }(t)&=&-\left[\frac{h^{\prime }(T-t)}{h(T-t)}+K\right]%
a(t)+ (\gamma-1)(a(t))^{\frac{\gamma}{\gamma-1}} \\
&+&\int_{t}^{T} h(T-t)\frac{\partial}{\partial t}\left[\frac{h(s-t)}{h(T-t)}%
\right](a(s))^{\frac{\gamma}{\gamma-1}}e^{-\left(\int_{t}^{s}\gamma(a(u))^{%
\frac{1}{\gamma-1}}\,du\right)}\,ds.  \notag
\end{eqnarray}

We notice that on $[T-\epsilon, T]$

\begin{equation*}
a^{\prime }(t)\approx -\left[\frac{h^{\prime }(T-t)}{h(T-t)}+K\right]%
a(t)+ (\gamma-1)(a(t))^{\frac{\gamma}{\gamma-1}} +O(\epsilon).
\end{equation*}

Consequently, for the choice of our parameters in the numerical experiment
we get that $\frac{h^{\prime }(T-t)}{h(T-t)}+K<-1$ for small $\epsilon$
(keeping in mind that $a(1)=2$) we see that $a(t)$ is increasing
on $[T-\epsilon, T].$ For hyperbolic discounting it can be shown that
$\frac{\partial}{\partial t}\left[\frac{h(s-t)}{h(T-t)}\right]<0.$
It is obvious that $a(t)$ is decreasing
in a neighborhood of 0 due to the negative
contribution of the term $$\int_{t}^{T} h(T-t)\frac{\partial}{\partial t}\left[\frac{%
h(s-t)}{h(T-t)}\right](a(s))^{\frac{\gamma}{\gamma-1}}e^{-\left(\int_{t}^{s}\gamma(%
a(u))^{\frac{1}{\gamma-1}}\,du\right)}\,ds.$$
In conclusion, the consumption rate policy, $\frac{F_{2}(t,x)}{x}=[a(t)]^{\frac{1}{\gamma -1}},$ (see \eqref{0109con}) is neither increasing nor decreasing in time.

\begin{flushright}
$\square$
\end{flushright}

{\bf {Appendix E}}: Proof of Lemma \ref{QuadraticSol}: Let $x\triangleq a^{\frac{1}{1-\g}}$, $\alpha_j\triangleq \lam+r_j-K-\frac{\g\eta}{l},$ $j=1,2.$ Equation \eqref{eq_a} becomes

\begin{equation}
\frac{1}{x}=\frac{1}{\aone+\g\beta x }+\frac{\lam m^{\frac{1}{1-\gamma}}}{\atwo+\g\beta x }
\end{equation}

For the sake of simplicity, consider the case $m=1$.

 We want to find $x>0$ which solves
 \begin{equation}\label{Q}
\g \beta (1-\g\beta+\lam) x^2+ (\atwo(1-\g\beta)+\aone(\lam-\g\beta) ) x- \aone\atwo=0
\end{equation}
 and  transversality conditions \eqref{TC} i.e
 \begin{equation}\label{TC2}
 \aone+\g\beta x>0\\\qquad
 \atwo+\g\beta x>0
 \end{equation}

 This can be done by splitting the analysis into 3 cases: $\gamma\in(0,1), \gamma=0,$ and $\gamma<0.$ We omit the details.

\begin{flushright}
$\square$
\end{flushright}

{\bf {Appendix F}}: Proof of Lemma \ref{Le1}: Let $a=inf_{t\in\left[0,T\right]} a(t)$. Lemma \ref{existenceODE} guarantees that $a>0.$
We show that

 \begin{equation}\label{120}
 a_n^1\geq \frac{2A}{3}\quad \forall n\in  0,1,\cdots, N,
 \end{equation}

 and this makes the recursive scheme well defined. We prove \eqref{120} by mathematical induction. Assume that

 \begin{equation}\label{112}
 a_k^1\geq \frac{2A}{3}\quad \forall k\in 0,1,\cdots, n,
 \end{equation}

 and prove that $a_{n+1}^1\geq \frac{2a}{3}.$ Let us define

 $$e_n \triangleq a_n^1 - a(t_n),\quad  f_n \triangleq A_n^1 - A(t_n)$$.

 By considering a second order Taylor expansion of $a(t_{n+1})$ at  $a(t_{n})$, we get $$a(t_{n+1}) = a(t_n) + \e a'(t_n) + c_n\e^2$$ with
$c_n$ a constant depending on $a$ and bounded by $c$ independently of $n$. Consequently
 \begin{eqnarray*}
 e_{n+1} &=& a_{n+1}^1 - a(t_{n+1})\\
 &=& a_n^1+\e\big((\g M(t_n)-\lambda(t_n)-1)(a_n^1)^{\frac{\g}{\g-1}}\\
&+&\left(\lambda(t_n)-\frac{h'(T-t_n)}{h(T-t_n)}-K\right)a_n^1+ \int_{t_n}^{T}L(s,t_n)(a(s))^{\frac{\g}{\g-1}}\left(\frac{A(s)}{A(t_n)}\right)ds\\
&-&\bigg(a(t_n)+ \e (\g M(t_n)-\lambda(t_n)-1)(a(t_n))^{\frac{\g}{\g-1}}\\
&+&\left(\lambda(t_n)-\frac{h'(T-t_n)}{h(T-t_n)}-K\right)a(t_n)\\
&+& \int_{t_n}^{T}L(s,t_n)(a(s))^{\frac{\g}{\g-1}}\left( \frac{A(s)}{A(t_n)}\right)ds+ c_n\e^2\bigg)\\
&=&e_n +\e (\g M(t_n)-\lambda(t_n)-1)(  (a_n^1)^{\frac{\g}{\g-1}}-(a(t_n))^{\frac{\g}{\g-1}}  )\\
&+&\e \left(\lambda(t_n)-\frac{h'(T-t_n)}{h(T-t_n)}-K\right)e_n -c_n \e^2
 \end{eqnarray*}

By the mean value Theorem applied to the function $x\rightarrow x^{\frac{\g}{\g-1}}, $ one gets that

 $$|(a_k^1)^{\frac{\g}{\g-1}}-(a(t_k))^{\frac{\g}{\g-1}}|\leq {\frac{\g}{\g-1}}\left(\frac{2a}{3}\right)^{\frac{1}{\g-1}} |a_k^1-a(t_k)|.$$

 Therefore there exists $M>0$, such that
 \begin{equation}\label{II}
 |e_{k+1}|\leq |e_k|(1+M|\e|)+c\e^2,\quad \forall k\in0,1,\cdots, n.
 \end{equation}

 By iterating \eqref{II} for $k=0\cdots n,$ one gets

 $$|e_{n+1}|\leq c\e^2 \frac{(1+M \e)^{n+1}-1}{M\e}\leq|c|\e^2 \frac{e^{MT}-1}{\frac{MT}{N}}\leq C |\e|,$$

 for some constant $C$ independent of $n.$ Therefore

 $$a^1_{n+1} \geq a(t_{n+1})-|e_{n+1}|\geq a-C|\e|\geq 2a/3 $$ for $|\e|$ small enough. This proves \eqref{120}.

 Moreover it follows that

  $$|a_n^1-a(t_n)|=|e_n|\leq C|\e|, \; \forall n\in 0,1,\cdots,N.$$

Similar arguments show that

$$|A_n^1-A(t_n)| \leq C|\e|,\quad \forall n\in   0,1,\cdots, N.$$

\begin{flushright}
$\square$
\end{flushright}

{\bf {Appendix G}}: Proof of Lemma \ref{Le2}:  We show that

 \begin{equation}\label{0120}
 a_n^2\geq \frac{a}{2}\quad \forall n\in  0,1,\cdots, N,
 \end{equation}

 and this makes the recursive scheme well defined. We prove \eqref{0120} by mathematical induction.
  Assume that
 \begin{equation}\label{0112}
 a_k^2\geq \frac{a}{2}\quad \forall k\in 0,1,\cdots, n,
 \end{equation}

 and prove that $a_{n+1}^2\geq \frac{A}{2}.$ Let $r_n\triangleq a_n^2-a_n^1,$ so
  \begin{eqnarray*}
 r_{n+1} &=& a_{n+1}^2-a_{n+1}^1\\
 &=&a_n^2+\e(\g M(t_n)-\lambda(t_n)-1)(a_n^2)^{\frac{\g}{\g-1}}+\e\left(\lambda(t_n)-\frac{h'(T-t_n)}{h(T-t_n)}-K\right)a_n^2\\
&-&\e^2 \sum_{j=0}^{n-1}L(t_j,t_n)(a(t_j))^{\frac{\g}{\g-1}}A(t_j)-\e \int_{t_n}^{T}L(s,t_n)(a(s))^{\frac{\g}{\g-1}}\left(\frac{A(s)}{A(t_n)}\right)ds \\
 &-& a_n^1-\e(\g M(t_n)-\lambda(t_n)-1)(a_n^1)^{\frac{\g}{\g-1}}
-\e\left(\lambda(t_n)-\frac{h'(T-t_n)}{h(T-t_n)}-K\right)a_n^1\\
 &=&r_n+\e(\g M(t_n)-\lambda(t_n)-1)\big( (a_n^2)^{\frac{\g}{\g-1}}-(a_n^1)^{\frac{\g}{\g-1}}\big)\\
 &+&\e\left(\lambda(t_n)-\frac{h'(T-t_n)}{h(T-t_n)}-K\right)r_n \\
 &+&\e \sum_{j=0}^{n-1}\left[-\e L(t_j,t_n)(a(t_j))^{\frac{\g}{\g-1}}\left(\frac{A(t_j)}{A(t_n)}\right)-\int_{t_{j+1}}^{t_{j}}L(s,t_n)(a(s))^{\frac{\g}{\g-1}}\left(\frac{A(s)}{A(t_n)}\right)ds\right] \\
 \end{eqnarray*}

 Moreover
 \begin{eqnarray*}
 &&\big|-\e L(t_j,t_n)(a(t_j))^{\frac{\g}{\g-1}}\left(\frac{A(t_j)}{A(t_n)}\right)-\int_{t_{j+1}}^{t_{j}}L(s,t_n)(a(s))^{\frac{\g}{\g-1}}\left(\frac{A(s)}{A(t_n)}\right)ds\big|\\
 &=&\big|\int_{t_{j+1}}^{t_{j}}\left(L(t_j,t_n)(a(t_j))^{\frac{\g}{\g-1}}\left(\frac{A(t_j)}{A(t_n)}\right)-L(s,t_n)(a(s))^{\frac{\g}{\g-1}}\left(\frac{A(s)}{A(t_n)}\right)\right) ds\big|\\
 &\leq&\frac{1}{A(t_n)}\big|\int_{t_{j+1}}^{t_{j}} \bigg( [L(t_j,t_n)(a(t_j))^{\frac{\g}{\g-1}}A(t_j)-L(s,t_n)(a(t_j))^{\frac{\g}{\g-1}}A(t_j)]+[L(s,t_n)(a(t_j))^{\frac{\g}{\g-1}}A(t_j)\\
 &-&L(s,t_n)(a(s))^{\frac{\g}{\g-1}}A(t_j)]+[L(s,t_n)(a(s))^{\frac{\g}{\g-1}}A(t_j)-L(s,t_n)(a(s))^{\frac{\g}{\g-1}}A(s)]\bigg)ds\big| \\
 &\leq&\frac{1}{A(t_n)}\big( K_0(t_j-t_{j+1})^2+K_1(t_j-t_{j+1})^2+K_2(t_j-t_{j+1})^2\big)\leq K_4 \e^2,
\end{eqnarray*}

for some positive constants $K_0, K_1, K_2, K_4.$ The last inequalities follow from the boundedness of $a(t)$ ( see Lemma \ref{existenceODE})
 and from the boundedness  of coefficients. Arguing as in Lemma \ref{Le1}, we can then find $M>0$ such that

  \begin{equation}\label{1II}
 |r_{k+1}|\leq |r_k|(1+M|\e|)+c\e^2,\quad \forall k\in0,1,\cdots, n.
 \end{equation}

 By iterating \eqref{1II} for $ k=0\cdots n,$ one gets

 $$|r_{n+1}|\leq c\e^2 \frac{(1+M \e)^{n+1}-1}{M\e}\leq|c|\e^2 \frac{e^{MT}-1}{\frac{MT}{N}}\leq C |\e|,$$

 for some constant $C$ independent of $n.$ Therefore

 $$a^2_{n+1} \geq a^1_{n+1}-|r_{n+1}|\geq 2a/3-C|\e|\geq a/2 $$ for $|\e|$ small enough. This proves \eqref{0120}.  Moreover it follows that

  $$|a_n^2-a_n^1|=|r_n|\leq C|\e|, \; \forall n\in 0,1,\cdots,N.$$

Similar arguments show that

$$|A_n^2-A_n^1| \leq C|\e|,\quad \forall n\in   0,1,\cdots, N.$$

\begin{flushright}
$\square$
\end{flushright}

{\bf {Appendix H}}: Proof of Lemma \ref{Le3}: We show that

 \begin{equation}\label{00120}
 a_n^3\geq \frac{a}{4}\quad \forall n\in  0,1,\cdots, N,
 \end{equation}

 and this makes the recursive scheme well defined. We prove \eqref{00120} by mathematical induction.
  Assume that
 \begin{equation}\label{00112}
 a_k^3\geq \frac{a}{4}\quad \forall k\in 0,1,\cdots, n,
 \end{equation}

 and prove that $a_{n+1}^3\geq \frac{a}{4}.$ Let us introduce $u_n\triangleq a_n^3-a_n^2$ and $v_n\triangleq A_n^3-A_n^2$. It follows that
  \begin{eqnarray*}
 u_{n+1} &=& a_{n+1}^3-a_{n+1}^2\\
 &=&a_n^3+\e(\g M(t_n)-\lambda(t_n)-1)(a_n^3)^{\frac{\g}{\g-1}}+\e\left(\lambda(t_n)-\frac{h'(T-t_n)}{h(T-t_n)}-K\right)a_n^3\\
&-&\e^2 \sum_{j=0}^{n-1}L(t_j,t_n)(a_j^3)^{\frac{\g}{\g-1}}A_j^3-a_n^2-\e(\g m(t_n)-\lambda(t_n)-1)(a_n^2)^{\frac{\g}{\g-1}}\\
&-&\e\left(\lambda(t_n)-\frac{h'(T-t_n)}{h(T-t_n)}-K\right)a_n^2+\e^2 \sum_{j=0}^{n-1}L(t_j,t_n)(a(t_j))^{\frac{\g}{\g-1}}A(t_j)\\
&=&u_n+\e(\g m(t_n)-\lambda(t_n)-1)( (a_n^3)^{\frac{\g}{\g-1}}- (a_n^2)^{\frac{\g}{\g-1}})\\
&+& \e\left(\lambda(t_n)-\frac{h'(T-t_n)}{h(T-t_n)}-K\right)u_n^3-\e^2 \sum_{j=0}^{n-1}L(t_j,t_n)r_{j,n},
 \end{eqnarray*}

 where $r_{j,n}\triangleq (a_j^3)^{\frac{\g}{\g-1}}\left(\frac{A_j^3}{A_n^3}\right)- (a(t_j))^{\frac{\g}{\g-1}}\left(\frac{A(t_j)}{A(t_n)}\right).$ By triangle inequality it follows that
  \begin{eqnarray*}
  |r_{j,n}|&\leq&|(a_j^3)^{\frac{\g}{\g-1}}\left(\frac{A_j^3}{A_n^3}\right)-(a_j^2)^{\frac{\g}{\g-1}}\left(\frac{A_j^3}{A_n^3}\right)|+|(a_j^2)^{\frac{\g}{\g-1}}\left(\frac{A_j^3}{A_n^3}\right)- (a_j^2)^{\frac{\g}{\g-1}}\left(\frac{A_j^2}{A_n^3}\right)|\\
  &+&|(a_j^2)^{\frac{\g}{\g-1}}\left(\frac{A_j^2}{A_n^3}\right)-(a_j^2)^{\frac{\g}{\g-1}}\left(\frac{A_j^2}{A_n^2}\right)|+|(a_j^2)^{\frac{\g}{\g-1}}\left(\frac{A_j^2}{A_n^2}\right)-(a(t_j))^{\frac{\g}{\g-1}}\left(\frac{A(t_j)}{A(t_n)}\right)|\\
  &\leq&M_1|u_j|+M_2|v_j|+M_3|v_n|+M_4|\e|,
 \end{eqnarray*}

 for some positive constants $M_1, M_2, M_3, M_4,$ where the last inequality from the boundedness of $a(t)$ ( see Lemma \ref{existenceODE})
 and of other coefficients in our model. Arguing as in the previous Lemmas one can find a positive constant $C$ such that
 \begin{equation}\label{M}
 |u_{n+1}|\leq |u_n|+C|\e u_n|+C|\e|(max_{j\in 0,\cdots,n }|u_j|+max_{j\in 0,\cdots,n}|v_j|)+C\e^2.
 \end{equation}

 On the other hand

 $$v_{n+1}=v_n -\g \e M(t_n) \big((a_n^3)^{\frac{\g}{\g-1}}p_n^3-(a_n^2)^{\frac{\g}{\g-1}}p_n^2\big),$$

 hence

  $$|v_{n+1}|\leq|v_n|+|\g\e M(t_n)|\big( (a_n^3)^{\frac{\g}{\g-1}}|A_n^3-A_n^2| + A_n^2 |(a_n^3)^{\frac{\g}{\g-1}}-(a_n^2)^{\frac{\g}{\g-1}}| \big).$$

  However this implies that
 \begin{equation}\label{M1}
 |v_{n+1}|\leq |v_n|+M|\e|(|u_n|+|v_n|),
 \end{equation}

 for some positive constant $M.$ Let us define $x_n=max_{j\in 0,\cdots,n }|u_j|$ and $y_n=max_{j\in 0,\cdots,n}|v_j| $
and $z_n=x_n + y_n$. Inequalities \eqref{M} and \eqref{M1} hold also for $k\leq n,$ i.e.,

$$ |u_{k+1}|\leq |u_k|+C|\e u_k|+C|\e| (x_k + y_k)+C\e^2,\,\,|v_{k+1}|\leq |v_k|+M|\e|(|u_k|+|v_k|).$$

By taking maximum over $k\in 0,\cdots,n$ in these inequalities one obtains

$$x_{n+1}\leq x_n+2C|\e| x_n+C|\e|y_n+C\e^2$$ and

 $$y_{n+1}\leq y_n + M|\e| (x_n+y_n)$$

By adding these inequalities, it follows that

 $$z_{n+1}\leq z_n+(2C+M)|\e|z_n+M\e^2.$$

This in turn yields that $z_n\leq C |\e|,$ for some positive constant still denoted (with some abuse of notations) $C.$  Therefore

 $$a^3_{n+1} \geq a^2_{n+1}-|u_{n+1}|\geq a/2-K|\e|\geq a/4 $$ for $|\e|$ small enough. This proves \eqref{00120}.  Moreover  $z_n\leq C |\e|,$

 implies that

  $$|a_n^3-a_n^2|=|u_n|\leq C|\e|, \; \forall n\in 0,1,\cdots,N.$$

and

$$|A_n^3-A_n^2|=|v_n| \leq C|\e|,\quad \forall n\in   0,1,\cdots, N.$$

\begin{flushright}
$\square$
\end{flushright}

\end{document}